\documentclass[12pt]{amsart}
\pdfoutput=1

\usepackage{amscd,latexsym,amsthm,amsfonts,amssymb,amsmath,amsxtra,}
\usepackage[mathscr]{eucal}
\renewcommand\theequation{\thesection.\arabic{equation}}

\newcommand{\BA}{{\mathbb {A}}}

\newcommand{\BQ}{{\mathbb {Q}}}
\newcommand{\BR}{{\mathbb {R}}}

\newcommand{\BZ}{{\mathbb {Z}}}

\newcommand{\CA}{{\mathcal {A}}}

\newcommand{\CO}{{\mathcal {O}}}

\newcommand{\GL}{{\mathrm{GL}}}

\renewcommand{\Im}{{\mathrm{Im}}}

\newcommand{\PGL}{{\mathrm{PGL}}}

\renewcommand{\Re}{{\mathrm{Re}}}

\newcommand{\SL}{{\mathrm{SL}}}

\newcommand{\bs}{\backslash}

\def\eps{{\epsilon}}

\newtheorem{thm}{Theorem}[section]

\newtheorem{lem}[thm]{Lemma}
\newtheorem{prop}[thm]{Proposition}

\newtheorem {ques/conj}[thm]{Question/Conjecture}

\newtheorem{defn}[thm]{Definition}
\newtheorem{rmk}[thm]{Remark}

\newtheorem{exmp}[thm]{Example}

\makeatletter

\newcommand{\Rmnum}[1]{\expandafter\@slowromancap\romannumeral #1@}
\makeatother

\begin{document}
\renewcommand{\theequation}{\arabic{equation}}
\numberwithin{equation}{section}

\title{The Weyl bound for triple product $L$-functions in the cubic level}

\author{Xinchen Miao}
\address{Mathematisches Institut\\ Endenicher Allee 60, Bonn, 53115, Germany}
\address{Max Planck Institute for Mathematics\\ Vivatsgasse 7, Bonn, 53111, Germany}
\email{miao@math.uni-bonn.de, olivermiaoxinchen@gmail.com}

\author{Huimin Zhang}
\address{School of Mathematical Sciences, Zhejiang University \\
Hangzhou, Zhejiang 310058, China}
\email{huimin.zhang@zju.edu.cn}

\date{August, 2025}

\begin{abstract}

 In this paper, we focus on the strong subconvexity bounds for triple product $L$-functions in the cubic level aspect. 
Our proof on the Weyl-type bound synthesizes techniques from classical analytic number theory with methods in automorphic forms and representation theory.
The methods include the refined Petersson trace formula for the newforms of cubic level,
classical Voronoi summation formula, Jutila's circle method, Kuznetsov trace formula and the spectral large sieve inequality. 

\end{abstract}

\thanks{The first author is supported by ERC Advanced Grant  101054336 and Germany's Excellence Strategy grant EXC-2047/1 - 390685813.}

\thanks{The second author is in part supported by the National Key R$\&$D Program of China (No. 2021YFA1000700), NSFC (No. 12031008) and China Scholarship Council (No. 202206220071).}
	
\keywords{triple product $L$-function, Weyl bound, cubic level, the refined Petersson trace formula, shifted convolution problem}
	
\subjclass[2010]{11F72, 11F70.}

\maketitle

\tableofcontents

\section{Introduction}

\subsection{Background and main theorem}

Subconvexity estimation is one of the most important and challenging problem in the theory of analytic number theory and $L$-functions. 
In general, let $F$ be a number field with adele ring $\BA_F$, and let $\Pi$ be an automorphic representation of a reductive group $\mathrm{G}$. Let $L(s, \Pi)$ be the corresponding automorphic $L$-function associated to $\Pi$. If $C(\Pi)$ denotes the analytic conductor of $L(s,\Pi)$, then the classical Phragm\'en--Lindel\"of principle gives the upper bound $C(\Pi)^{\frac{1}{4}+\epsilon}$ on the critical line $\Re(s)=\frac{1}{2}$. The subconvexity problem for $L(\frac{1}{2}, \Pi)$ is to establish a non-trivial upper bound of the following shape:
$$L(\frac{1}{2}, \Pi) \ll_{F,\epsilon} C(\Pi)^{\frac{1}{4}-\delta+\epsilon},$$
where $\delta$ is some positive absolute constant which is independent of $C(\Pi)$. If one can obtain $\delta = 1/12$, then we refer to as a Weyl-type subconvex bound.

In the lower-rank case $G=\GL_1, \GL_2$, the subconvexity problem has now been solved completely over a fixed general number field $F$, uniformly in all aspects ($t$, weight, spectral, level) \cite{subconvexity}. The main ingredients of the proof are integral representations, period integrals of certain $L$-functions (Ichino--formula), and local period integral computations.

In the higher-rank case, for example $G=\GL_2 \times \GL_2$ or $G=\GL_2 \times \GL_2 \times \GL_2$, the situation is far from well-understood (see \cite{michel2} \cite{Hu17} \cite{michel0} \cite{michel1}  \cite{sparse}). 
For more history on the subconvexity problem, the interested readers may see \cite{mic} for a survey and more details.

In this paper, we focus on the level aspect subconvexity problem of the triple product $L$-function $L(1/2, \pi_1 \times \pi_2 \times \pi_3)$ in the cubic level (horizontal) aspect, which is the case $G=\GL_2 \times \GL_2 \times \GL_2$. By local triple product integral computation and the amplification method, it is known that $L(1/2, \pi_1 \times \pi_2 \times \pi_3) \ll_{\varepsilon} p^{3-0.25+\varepsilon}$ in \cite{Hu17}. We improve the subconvexity bound and achieve the Weyl-type bound assuming the Ramanujan--Petersson conjecture (see Theorem \ref{triple}) in the cubic level case. Our method is a combination of different tools in classical analytic number theory, automorphic forms and representation theory, including Petersson trace formula for the newforms over a short family for cubic level, Voronoi summation formula, period integrals, Jutila's circle method, Kuznetsov trace formula and large sieve inequality.

We fix the ground number field $F=\BQ$. Let $\pi_1, \pi_2, \pi_3$ be three cuspidal automorphic representations of $\PGL_2(\BA_\BQ)$. We further assume that their corresponding weights are $k,k,2k$, where $k$ is a fixed even number. This means that the local archimedean component of the global representation is a discrete series representation of the given weight. We further assume that their corresponding levels are $1,1,p^3$, where $p>2$ is a prime number. For $\pi_1,\pi_2,\pi_3$ given above, we can realize them as three holomorphic cusp newforms $f$, $g$ and $h$ respectively. Hence, $f$ and $g$ are both holomorphic cusp forms with weight $k$ and level $1$, while $h$ is a holomorphic cusp form with weight $2k$ and level $p^3$. 
We will consider the subconvexity bound for the triple product $L$-function $L(1/2,f \times g \times h)$ when $p \rightarrow +\infty$.

For the local place $v=p$, the representation $\pi_{3,v}$ is a simple supercuspidal representation for $\PGL_2(\BQ_v)$ (\cite{KL15}). For the cubic level case, the Petersson or Kuznetsov trace formulae we apply will be the simplest (\cite{PWZ20}).

We let the real number $\theta$ be the best exponent toward the Ramanujan--Petersson Conjecture for $\GL(2)$ over $\BQ$. Hence, we have $0 \leq  \theta \leq \frac{7}{64}$. 

With the above notation, we can state our main result now.

\begin{thm} \label{triple}
For fixed even number $k>k_{\varepsilon}$, we have
$$ L\left(\frac{1}{2}, f \times g \times h \right) \ll_{k, \varepsilon} p^{2+2\theta+\varepsilon}.$$
\end{thm}

\begin{rmk}
Here $k_{\varepsilon}$ is an absolute constant that only depends on $\varepsilon$. We do not try to optimize the value for $k_{\varepsilon}$. It is possible to simply pick $k_{\varepsilon}=10^5/\varepsilon$ for any fixed $\varepsilon>0$. If we assume the Ramanujan--Petersson conjecture, i.e. $\theta=0$, for sufficiently large weight $k$, we roughly get the Weyl-type bound for the triple product $L$-function by noting that the convexity bound is $p^{3+\varepsilon}$. This is the first example to achieve the strong Weyl-type bound for higher degree automorphic $L$-functions in level aspect. This also improves the corresponding result in \cite{Hu17}, which proves the subconvexity bound $p^{3-0.25+\varepsilon}$. We note that $2+2\theta < 2+\frac{2}{9}<3 \times \frac{3}{4}$. Hence, our subconvexity bound is stronger than the Burgess bound unconditionally.
\end{rmk}

\begin{rmk}
We may also consider the general case in which the conductor for $\pi_3$ is $p^n$ for $n \geq 2$. In \cite{Hu17}, the subconvexity bound gives that $L(1/2, \pi_1 \times \pi_2 \times \pi_3) \ll_{\pi_1,\pi_2,\pi_{3,\infty}} p^{(11/12+\varepsilon)n}$ for all $p,n \geq 2$. It is possible to generalize this result by the application of the refined Petersson and Kuznetsov trace formula developed in \cite{Hu24}.
\end{rmk}

\subsection{Outline of methods}

Our methods include different and powerful tools in analytic number theory, automorphic forms and representation theory. Because of the lack of approximate functional equations, it is really hard to express the central value of triple product $L$-function $L(1/2,f \times g \times h)$ in terms of normalized Fourier coefficients $\lambda_f(n),\lambda_g(n)$ and $\lambda_h(n)$. However, for the special case that all the three cusp forms are holomorphic, we can overcome this obstacle by applying the theory of Poincar\'e series. It is well known that the space of holomorphic cusp forms with weight $k$ and level one can be generated by the first $m$ Poincar\'e series, where $m$ is the dimension of $S_k(1)$. Hence, without loss of generality, we can assume that $g$ is a level one Poincar\'e series with shift one. Scaling by $p^3$ which is the classical way to produce an old form suggests us to substitute the initial Poincar\'e series with a new Poincar\'e series with level $p^3$ and shift $p^3$. With this new Poincar\'e series, we can unfold the period integral and roughly achieve that the central value $L(1/2, f \times g \times h)$ equals the square of convolution sums in terms of $\lambda_f(n)$ and $\lambda_h(m)$ with the shift $m-n=p^3$. This finishes the setup part.

After we open the square, each term looks like $\lambda_f(n)\lambda_f(m)\lambda_h(n+p^3)\lambda_h(m+p^3)$. From the non-negativity property, we may consider $\sum_{h \in \CA(2k,M,p^3)} \lambda_h(n+p^3)\lambda_h(m+p^3)$. Since $\vert \CA(2k,M,p^3) \rvert \asymp kp^2$, which is a short family of $\CA(2k,p^3)$ with $\vert \CA(2k,p^3) \rvert \asymp kp^3$, this encourages us to apply the refined Petersson trace formula with the prescribed local component over a short family at the finite place $v=p$.

After the application of the refined Petersson trace formula which is a combination of automorphic forms and representation theory, we reach the summation with Hecke eigenvalues, Kloosterman sums and classical Bessel functions involved (see summation $\mathcal{O}_3$ for example in Section 3.1). Now the powerful tools in analytic number theory and automorphic forms play the most important role on the stage. We apply Voronoi summation formula in one variable. The resulting Kloosterman sums simplifies enormously because the shift $p^3$ equals the level. We apply Kuznetsov trace formula backwards and arrive at a shifted convolution problem that we can treat by Jutila's circle method along with Voronoi summation formula on both variables, Kuznetsov trace formula again and finally the spectral large sieve inequality, which will achieve the Weyl bound for triple product $L$-functions in the cubic level. This is the first example to achieve the strong Weyl bound for higher degree automorphic $L$-functions in the horizontal level aspect. For spectral aspect analogy, the Weyl bound for triple product $L$-functions was proven in \cite{BJN23}.

{\textbf{Acknowledgments}}
 The first author would like to express sincere gratitude to Professor Valentin Blomer and Professor Paul Nelson for stimulating discussions and for generously sharing their personal notes. The first author would like to thank the University of Bonn and the Max Planck Institute for Mathematics in Bonn for excellent working conditions and research environment. The author is supported by ERC Advanced Grant  101054336 and Germany's Excellence Strategy grant EXC-2047/1 - 390685813.
 The second author would like to express sincere gratitude to Chenchen Shao for helpful discussions.
The second author is deeply grateful to the University of Bonn and Shandong University
for providing a perfect research atmosphere, where parts of this article were worked out.  
And She is grateful to the China Scholarship Council for supporting her studies at the University of Bonn. The second author is in part supported by the National Key R$\&$D Program of China (No. 2021YFA1000700), NSFC (No. 12031008) and China Scholarship Council (No. 202206220071).

\medskip
{\textbf{Notation}}
Throughout the paper, $A$ denotes arbitrarily large positive constants and $\varepsilon$ are sufficiently small positive constants, not necessarily the same at different occurrences. 
We use $A\asymp B$ to mean that $c_1B\leq |A|\leq c_2B$ for some positive constants $c_1$ and $c_2$. The symbol $\ll_{a,b,c}$ denotes that the implied constant depends at most on $a$, $b$ and $c$. In addition, $V$ denotes generally smooth weight
functions, not necessarily the same at different occurrences. And $S(m,n,c)$ denotes the classical Kloosterman sums.
As usual, $e(x)=e^{2 \pi i x}$ and
$q\sim R$ means $R \leq q < 2R$.

\section{Preliminaries and Background}

\subsection{Automorphic forms and $L$-functions}\label{autl} 
Let $S_k(M)$ denote the space of holomorphic cusp forms of level $M$, weight $k$, and trivial nebentypus. Let $S_k^*(M)\subset S_k(M) $ denote the space of newforms. Every $f\in S_k(M)$ has a Fourier series expansion
\begin{align*}
f(z)=\sum_{n=1}^\infty \lambda_f(n) n^{\frac{k-1}{2}} e(nz)
\end{align*}
for $\Im(z)>0$. Let $B_k(M)$ denote an orthogonal basis of $S_k(M)$, which contains a basis $B_k^*(M)$ of $S_k^*(M)$, normalized so that $\lambda_f(1)=1$ for every $f\in B_k^*(M)$. Thus for $f\in B_k^*(M)$, the $n$-th Fourier coefficient equals the $n$-th Hecke eigenvalue.

We now consider more general automorphic forms. Let $N\ge 1$ be an integer, $\psi$ a Dirichlet character of modulus $N$, and let cond$(\psi$) denote the modulus of the primitive character which induces $\psi$. We further suppose that $\psi$ is even and cond$(\psi$) is squarefree and odd, since this is enough for our purpose. Actually, throughout our paper, we have cond$(\psi)=1$. Let $B_k(N,\psi)$ denote an orthogonal basis of holomorphic cusp forms of weight $k$ and nebentypus $\psi$ with respect to the standard congruence subgroup $\Gamma_0(N)$. Let $B(N,\psi)$ denote an orthogonal basis of Maass cusp forms with nebentypus $\psi$ with respect to $\Gamma_0(N)$, with $t_g$ denoting the spectral parameter of $g\in B(N,\psi)$. It is conjectured that $t_g$ is real (which is known as the Selberg eigenvalue conjecture), but the possibility of purely imaginary $t_g$ with $|t_g|\in(0,1/2)$ has not been disproven (these are called exceptional eigenvalues). If the level $1 \leq N \leq 857$, it is known that there does not exist exceptional eigenvalues. If we work with orthonormal basis, we will divide by the $L^2$-norm $\| g\|$ of the forms. Let $E_{\mathfrak{c}}(\cdot, 1/2+it)$ denote the Eisenstein series of nebentypus $\psi$ associated to a cusp $\mathfrak{c}$ which is singular for $\psi$. For $\sigma=
(\begin{smallmatrix}
* & * \\
c & d 
\end{smallmatrix})
\in SL_2(\mathbb{R})$, let $j(\sigma, z):=cz+d$. We write the Fourier expansions of these objects around a singular cusp $\mathfrak{a}$ associated to a scaling matrix $\sigma_\mathfrak{a}$ as follows:
\begin{align*}
&g(\sigma_\mathfrak{a} z)j(\sigma_\mathfrak{a}, z)^{-k}=\sum_{n=1}^{\infty} \lambda_{g_\mathfrak{a}}(n) n^{\frac{k-1}{2}} e(nz) \ \ \text{for } g\in B_k(N,\psi),\\
&g(\sigma_\mathfrak{a}z)= \sum_{n\neq 0} \lambda_{g_\mathfrak{a}}(n) n^{-\frac{1}{2}} W_{0,it_f} (4\pi |n| y)e(nx) \ \ \text{for } g\in B(N,\psi),\\
&E_\mathfrak{c}(\sigma_\mathfrak{a}z,1/2+it)=c_{1,\mathfrak{c}}(t)y^{1/2+it}+c_{2,\mathfrak{c}}(t)y^{1/2-it}+\sum_{n\neq 0} \lambda_{\mathfrak{c}\mathfrak{a}}(n,t) n^{-\frac{1}{2}} W_{0,it_f} (4\pi |n| y)e(nx),
\end{align*}
where 
\[
W_{0,it_f} (4\pi |n| y)= \sqrt{n \pi y}K_{it_f}(2\pi |n| y)
\]
is a Whittaker function and $K_{it}$ is a Bessel function. For $\mathfrak{a}=\infty$, we simply write 
\[
\lambda_{g_\mathfrak{a}}(n)=\lambda_{g}(n) \ \ \text{ and } \ \ \lambda_{\mathfrak{c}\mathfrak{a}}(n,t)=\lambda_{\mathfrak{c}}(n,t).
\]

Let $B_k^*(N,\psi)\subset B_k(N,\psi)$ and $B^*(N,\psi)\subset B(N,\psi)$ denote orthogonal basis of the space of newforms, normalized so that $\lambda_g(1)=1$. For $g$ in $B_k^*(N,\psi)$ or $B^*(N,\psi)$, we have the Hecke multiplicativity relations
\[
\lambda_g(n)\lambda_g(m) =\sum_{\substack{d|(n,m)}} \psi(d) \lambda_g\Big(\frac{nm}{d^2}\Big), \ \ \ \lambda_g(nm)=\sum_{d|(n,m)} \mu(d) \psi(d) \lambda_g\Big(\frac{n}{d}\Big)\lambda_g\Big(\frac{m}{d}\Big),
\]
where $\mu$ is the Mobius function. For $g$ in $B_k^*(N,\psi)$, we have the Ramanujan bound (due to Deligne),
\[
\lambda_g(n)  \ll_{\varepsilon} n^\varepsilon,
\]
where $d(n)$ is the divisor function and for $g\in B^*(N,\psi)$, we have the following Kim--Sarnak bound 
\[
\lambda_g(n) \ll_{\varepsilon} n^{\frac{7}{64}+\varepsilon}.
\]

It will be useful to have basis $B_k(N,\psi)$ and $B(N,\psi)$ which are given in terms of lifts of newforms.

\begin{lem} \cite[Lemma 2.1]{Khan22} \label{specialbasis} We have an orthogonal basis
\begin{align}
B_k(N,\psi)= \bigcup_{\substack{\delta \ell | N }} \{ g^{(\delta)}: g\in B_k^*(\ell, \psi) \},
\end{align}
where
\[
g^{(\delta)}(z)=\sum_{r|\delta} \nu_\delta(r) r^{1/2} g(rz)
\]
for some complex numbers $\nu_\delta(r)$ depending on $r, \delta$, and $g$, such that $ \nu_\delta(r)\ll (r\delta)^\varepsilon$, where the understanding is that $B_k^*(\ell, \psi)$ is empty if $\psi$ is not a character mod $\ell$. Furthermore, this basis is orthonormal if every $g\in B_k^*(\ell, \psi)$ for $\ell|N$ is $L^2$-normalized with respect to $\Gamma_0(N)$.

\end{lem}

\subsection{The refined Petersson trace formula for the short moment}
One of the important ingredients for Theorem \ref{triple} is the refined Petersson trace formula for a short family over cubic level in \cite[Proposition 3.1]{PWZ20} (see also \cite[Theorem 1.1]{PWZ20}. 

\begin{thm} \cite[Proposition 3.1]{PWZ20} \label{peter1}
 For $(n_1n_2,p)=1$, we have
 \begin{equation}
 \begin{aligned}
  \sum_{\pi \in \CA(2k,p^3,M)}  \frac{\lambda_{\pi}(n_1)\lambda_{\pi}(n_2)}{L_{\text{fin}}(1,\pi,\text{sym}^2)}&= \delta(n_1,n_2)\cdot \frac{(2k-1)p^2}{2\pi^2} \\
  &+\frac{(-1)^k(2k-1)}{\pi} \cdot \sum_{c \geq 1} \frac{A_{p,M}(c)}{c} \cdot S(n_1,n_2, p^2c) \cdot J_{2k-1} \left( \frac{4\pi \sqrt{n_1n_2}}{p^2c} \right),
 \end{aligned}
 \end{equation}
\end{thm}
where $A_{p,M}(c)=e\left( \frac{M}{p}\right)$ if $p \nmid c$ and $A_{p,M}(c)=1$ if $p \vert c$. If $p \vert n_1$ or $p \vert n_2$, we have $\lambda_{\pi}(n_1)=0$ or $\lambda_{\pi}(n_2)=0$. $\CA(2k,p^3)$ is the set of holomorphic cusp newforms with weight $2k$, level $p^3$ and trivial nebentypus. $\CA(2k,p^3,M):=\{ \pi \in \CA(2k,p^3) \vert \pi_p=\pi_{M,\,\zeta_p},\; \zeta_p  \in \{\pm 1 \}\}$, where $\zeta_p$ (the local root number) and $M$ ($1 \leq M \leq p-1$) parametrize the simple supercuspidal representations. Hence, there are totally $2(p-1)$ simple supercuspidal representations for $\PGL_2(\BQ_p)$ up to isomorphism. For more explanations on $\CA(2k,p^3,M)$ and simple supercuspidal representations, we refer to \cite{PWZ20} and \cite{KL15}. We note that $\CA(2k,p^3) \asymp kp^3$ and $\CA(2k,p^3,M) \asymp kp^2$ for fixed $M$. Moreover, when we vary $M$ satisfying $1 \leq M \leq p-1$, we see that $\CA(2k,p^3)$ is the disjoint union of $\CA(2k,p^3,M)$ (\cite{PWZ20}). For a generalization for Theorem \ref{peter1} to level $p^n$ with $n \geq 2$, readers may refer to \cite{Hu24}, \cite{hpy24} and \cite{hpy25}.

Taking the summation for $M$ from $1 \leq M \leq p-1$, we have

\begin{thm} \cite[Theorem 1.1]{PWZ20} \label{peter3}
 For $(n_1n_2,p)=1$, we have
 \begin{equation}
 \begin{aligned}
  \sum_{\pi \in \CA(2k,p^3)}  \frac{\lambda_{\pi}(n_1)\lambda_{\pi}(n_2)}{L_{\text{fin}}(1,\pi,\text{sym}^2)}&= \delta(n_1,n_2)\cdot \frac{(2k-1)p^2(p-1)}{2\pi^2} \\
  &+\frac{(-1)^k(2k-1)}{\pi} \cdot \sum_{c \geq 1} \frac{A_{p}(c)}{c} \cdot S(n_1,n_2, p^2c) \cdot J_{2k-1} \left( \frac{4\pi \sqrt{n_1n_2}}{p^2c} \right),
 \end{aligned}
 \end{equation}
\end{thm}
where $A_{p}(c)=-1$ if $p \nmid c$ and $A_{p}(c)=p-1$ if $p \vert c$. If $p \vert n_1$ or $p \vert n_2$, we will have $\lambda_{\pi}(n_1)=0$ or $\lambda_{\pi}(n_2)=0$.

\subsection{Jutila's circle method}
We quote Jutila's circle method \cite{Ju95}.
\begin{lem}\label{lem3}  Let $Q \geq 1$ and $V$ be a smooth, nonnegative, nonzero  function with support in $[1, 2]$. For $r \in \mathbb{Q}$ write $I_{r}(\alpha)$ for the characteristic function of the interval $[r- 1/Q, r+1/Q]$ and define
\begin{displaymath}
  \Lambda := \sum_q V\Big(\frac{q}{Q}\Big) \phi(q),  \quad I(\alpha) = \frac{Q}{2   \Lambda} \sum_q V\Big(\frac{q}{Q}\Big) 
  \sum_{\substack{d\, (\text{{\rm mod }} q)\\ (d, q) = 1}} I_{d/q}(\alpha).
\end{displaymath}
Then $I(\alpha)$ 
is a good approximation to the characteristic function on $[0, 1]$ in the sense that
\begin{displaymath}
  \int_0^1 (1 - I(\alpha))^2 d\alpha \ll_{\eps} Q^{\eps - 1}
\end{displaymath}
for any $\eps > 0$.
\end{lem}

\subsection{Voronoi summation formula}
We need the following simple version of Voronoi summation formula.

\begin{lem} \label{voronoi} \cite[Lemma 3.3]{KMS24}
Let $f \in S_k(1)$ be a holomorphic Hecke cusp form with Fourier coefficients $\lambda_f(n)$ and trivial nebentypus. Let $a$ and $q$ be coprime integers. Let $g$ be a compactly supported smooth bump function on $\BR$. Then we have
$$ \sum_{n=1}^{\infty} \lambda_f(n) e \left( \frac{an}{q} \right) g(n)= \frac{1}{q} \cdot \sum_{m=1}^{\infty} \lambda_f(m) e \left( \frac{-m \overline{a}}{q} \right) H \left( \frac{m}{q^2} \right),$$
where $a \overline{a} \equiv 1 \; (\,\text{mod} \, q)$, and
$H(y)= 2 \pi i^k \int_0^{\infty} g(x) J_{k-1} (4 \pi \sqrt{xy} ) dx,$
where $J_{k-1}$ is the $J$-Bessel function and $k$ is the weight of $f$.
\end{lem}

\subsection{Spectral large sieve}\label{specsieve}

We state the spectral large sieve inequality.
\begin{thm}\label{largesievethm} \cite[Theorem 2.3]{Khan22}
For any sequence $\{\alpha_n\}$ of complex numbers, we have
\begin{equation}
\begin{aligned}
&\sum_{\substack{k\le T\\ k\equiv 0\bmod 2}} \frac{(k-1)!}{(4\pi)^{k-1}}  \sum_{g\in B_k(N,1)}\frac{1}{\| g \|^2} \Big| \sum_{M\le m \le 2M} \alpha_m \lambda_{g}(m )\Big|^2 \ll \Big(T^2+ \frac{M^{1+\epsilon}}{N}\Big)\| \alpha \|^2,\\
& \sum_{\substack{g\in B(N,1)\\ \vert t_f \rvert \le T}} \frac{1}{\| g \|^2}  \frac{1}{\cosh(\pi t_g)} \Big| \sum_{M\le m \le 2M} \alpha_m \lambda_{g}(m)\Big|^2 \ll \Big(T^2+ \frac{M^{1+\epsilon}}{N}\Big)\| \alpha \|^2,\\
& \sum_{\mathfrak{c} \text{ singular} } \int_{-\infty}^\infty \frac{1}{\cosh(\pi t)} \Big| \sum_{M\le m \le 2M} \alpha_m \lambda_{\mathfrak{c}}(m,t)\Big|^2 dt \ll \Big(T^2+ \frac{M^{1+\epsilon}}{N}\Big)\| \alpha \|^2,
\end{aligned}
\end{equation}
where $\| g\|^2$ is the Petersson inner product of $g$ with respect to the congruence subgroup $\Gamma_0(N)$ and $\| \alpha \|^2=\sum_{M\le m \le 2M} |\alpha_m|^2 $.
\end{thm}

When the sums on the left-hand side of the large sieve inequality are restricted to newforms, normalized to have the first Fourier coefficient equal to $1$, we can use
\begin{align}
\label{obs} N^{1-\epsilon}k^{-\epsilon} \frac{(k-1)!}{(4\pi)^{k-1}} \ll \| g \|^2\ll N^{1+\epsilon}k^\epsilon \frac{(k-1)!}{(4\pi)^{k-1}}
\end{align}
for $g$ holomorphic, and 
\begin{equation}
N^{1-\epsilon}(1+|t_g|)^{-\epsilon} \frac{1}{\cosh(\pi t_g)} \ll \| g \|^2\ll N^{1+\epsilon}(1+|t_g|)^{\epsilon}\frac{1}{\cosh(\pi t_g)}
\end{equation}
for $g$ a Maass form. 

\subsection{Kuznetsov trace formula}
Follow the notation from the previous section.
Let $\phi$ be a smooth function compactly supported on the positive real line. Define the following integral transforms of $\phi$ with the kernel in terms of Bessel functions:
\begin{equation}
\begin{aligned}
&\dot{\phi}(k)=\frac{i^k}{\pi}\int_0^\infty J_{k-1}(x)\phi(x)\frac{dx}{x},\\
&\tilde{\phi}(t)=\frac{2\pi i}{\sinh(\pi t)}\int_0^\infty (J_{2it}(x)-J_{-2it}(x)) \phi(x) \frac{dx}{x},\\
&\check{\phi}(t)= 8\cosh(\pi t)\int_0^\infty K_{2it}(x) \phi(x) \frac{dx}{x}.
\end{aligned}
\end{equation}
We have the following estimates for these transforms.
\begin{lem} \cite[Lemma 2.1]{BHM07} \cite[Lemma 2.4]{Khan22} \label{transf-bounds}
(a) If $\phi(x)$ is supported on $0<X<x<2X$ and satisfies $\phi^{(j)}(x)\ll_j (X/Z)^{-j}$ for all $j\ge 0$, then for $t\in\mathbb{R}$ we have
\begin{align*}
\dot{\phi}(t), \tilde{\phi}(t), \check{\phi}(t)\ll \frac{1+|\log (X/Z)|}{1+X/Z}.
\end{align*}
For $\vert t \rvert \geq 1$, we have
\begin{equation*}
\dot{\phi}(t), \tilde{\phi}(t), \check{\phi}(t)\ll \left( \frac{Z}{\vert t \rvert} \right) \cdot \left( \frac{1}{\sqrt{\vert t \rvert}}+\frac{X}{\vert t \rvert} \right). 
\end{equation*}
For $\vert t \rvert \geq \max(2X,1)$, we have
\begin{equation*}
\dot{\phi}(t), \tilde{\phi}(t), \check{\phi}(t)\ll_k \left( \frac{Z}{\vert t \rvert} \right)^k \cdot \left( \frac{1}{\sqrt{\vert t \rvert}}+\frac{X}{\vert t \rvert} \right). 
\end{equation*}

(b) If $\phi(x)$ is supported on $0<X<x<2X$ and satisfies $\phi^{(j)}(x)\ll_j (\frac{X}{Z})^{-j}$ for all $j\ge 0$, then for $t\in(-\frac{i}{4}, \frac{i}{4})$, we have
\begin{align*}
\tilde{\phi}(t), \check{\phi}(t)\ll \frac{1+(\frac{X}{Z})^{-2|\Im t|}}{1+(\frac{X}{Z})}.
\end{align*}

(c) If $\phi(x)=e^{ax}\psi(x)$, where $a=\pm 1$, and $\psi(x)$ is supported on $1\le X<x<2X$ and satisfies $\phi^{(j)}(x)\ll_j X^{-j}$ for all $j\ge 0$, then for $t\in\mathbb{R}$, we have
\begin{align*}
\dot{\phi}(t),  \tilde{\phi}(t) \ll X^{-1/2+\varepsilon}.
\end{align*}
For $|t|> X^{1/2+\varepsilon}$, we have 
\[
\dot{\phi}(t), \tilde{\phi}(t) \ll_l (|t|+X)^{-l}X^\varepsilon
\]
for any $l\ge 0$.
\end{lem}

We define 
\begin{equation}
\begin{aligned}
&\mathcal{H}=\sum_{\substack{k>0\\ k\equiv 0 \bmod 2}} \sum_{g\in B_k(N,1)}  \dot{\phi}(k) \frac{(k-1)!}{\pi (4\pi)^{k-1}} \frac{\lambda_{g}(n)\overline{\lambda_{g}(m)}}{\| g\|^2}\\
&\mathcal{M}= \sum_{g\in B(N,1)} \frac{1}{\cosh(\pi t_g)} \tilde{\phi}(t_g) \frac{\lambda_{g}(n)\overline{\lambda_{g}(m)}}{\| g\|^2}, \ \ \ \\
&\mathcal{E}=\sum_{\mathfrak{c}} \int_{-\infty}^{\infty}  \frac{1}{4\pi \cosh(\pi t)} \tilde{\phi}(t) \lambda_{\mathfrak{c}}(n,t)\overline{\lambda_{\mathfrak{c}}(m,t)} dt\\
&\mathcal{M}^\prime=\sum_{g\in B(N,1)} \frac{1}{\cosh(\pi t_g)} \check{\phi}(t_g) \frac{\lambda_{g}(-n)\overline{\lambda_{g}(m)}}{\| g\|^2}\\
&\mathcal{E}^\prime=\sum_{\mathfrak{c}} \int_{-\infty}^{\infty}  \frac{1}{4\pi \cosh(\pi t)} \check{\phi}(t) \lambda_{\mathfrak{c}}(-n,t)\overline{\lambda_{\mathfrak{c}}(m,t)} dt.
\end{aligned}
\end{equation}
We now state a special case of Kuznetsov's theorem which we will need.
\begin{thm} \label{kuznetsov} Keep the notation of Section \ref{specsieve}. 
Let $\phi$ be a smooth function compactly supported on the positive real line. For $n,m\ge 1$ we have
\begin{align*}
&\sum_{c\geq 1} \frac{1}{cN}  S(m,n, cN)  \phi\Big(\frac{4\pi\sqrt{nm}}{cN}\Big)=\mathcal{H}+\mathcal{M}+\mathcal{E},\\
&\sum_{c\geq 1} \frac{1}{cN}  S(-m,n, cN) \phi\Big(\frac{4\pi\sqrt{nm}}{cN}\Big) =\mathcal{M}^\prime+\mathcal{E}^\prime.\\
\end{align*}
\end{thm}

\subsection{Inert functions and stationary phase}
\label{subsec:inertfunctions}

We mention some properties of inert functions in this section. Inert functions are special families of smooth functions characterised by certain derivative bounds. 
\begin{defn}
Let $\mathcal{F}$ be an index set. A family $\{w_T\}_{T \in \mathcal{F}}$ of smooth function supported on a product of dyadic intervals in $\mathbb{R}_{>0}^{d}$ is called $X$-inert if for each $j = (j_1,j_2,\cdots,j_d) \in \mathbb{Z}_{>0}^d$, we have
\begin{equation*}
    C(j_1,j_2,\cdots,j_d) := \sup_{T \in \mathcal{F}} \ \sup_{(x_1,x_2,\cdots,x_d) \in \mathbb{R}_{>0}^{d}}  X^{-j_1-\cdots -j_d} \left\vert x_1^{j_1}\cdots x_d^{j_d} {w_T}^{(j_1,\cdots,j_d)}(x_1,\cdots,x_d)\right\vert < \infty.
\end{equation*}
\end{defn}

We will often denote the sequence of constants $C(j_1,j_2,\dots,j_d)$ associated with this inert function as $C_{\mathcal{F}}$. 

We note that the requirements for the functions to be supported on dyadic intervals can be achieved by applying a partition of unity.

We give some examples to highlight how such families can be constructed and some basic properties.

\begin{exmp} [Dilation]
    \label{ex:InertFamily}
    Let $w(x_1,\cdots,x_d)$ be a fixed smooth function that is supported on $[1,2]^d$, and define,
    \begin{equation}
        w_{X_1,\cdots,X_d} (x_1,\cdots,x_d) = w\left({\frac{x_1}{X_1},\cdots,\frac{x_d}{X_d}} \right).
    \end{equation}
    Then with $\mathcal{F} = \{ T = (X_1,\cdots,X_d) \in \mathbb{R}_{>0}^d\}$, the family $\{ w_T\}_{T \in \mathcal{F}}$ is $1$-inert.
\end{exmp}

The following propositions can be checked easily:

\begin{exmp}[Oscillation]\label{WTexample2} With $w$ as in the previous example we let
\begin{equation}
\label{eq:WTexample2}
 W_{T}(x_1, \dots, x_d) = e^{i \lambda_1 x_1 + \dots + i \lambda_d x_d} w\Big(\frac{x_1}{X_1} , \ldots, \frac{x_d}{X_d} \Big),
\end{equation}
 but now $\mathcal{F} =\{T = (X_1, \dots, X_d, \lambda_1, \dots, \lambda_d)\}$.  It is easy to see $W_T$ is $X$-inert with $X=X_T = 1 + \max(|\lambda_1| X_1, \dots, |\lambda_d| X_d)$, 
but not $Y$-inert for any $Y=Y_T$ such that $Y_T/X_T \to 0$ as $X_T\to \infty$.
 \end{exmp}

\begin{exmp}[Products]
If $w$ is an $X$-inert function and $v$ is a $Y$-inert function, then their product $w\cdot v$ is a $\max{(X,Y)}$-inert function.
\end{exmp}

Suppose that $w_T(x_1,\cdots,x_d)$ is $X$-inert and is supported on $x_i \asymp X_i$. Let
\begin{equation}
    \label{eq:InertFourierdefn}
    \widehat{w}_T(t_1,x_2,\cdots,x_d) = \int_{-\infty}^{\infty} w_T(x_1,\cdots,x_d) e(-x_1t_1)dx_1,
\end{equation}
denote its Fourier transform in the $x_1$-variable.

We state the following result (see \cite[Proposition 2.6]{KPY19}) regarding the Fourier transform of inert functions.
\begin{prop}
\label{prop:FourierInert}
Suppose that $\{ w_T : T \in \mathcal{F}\}$ is a family of $X$-inert functions such that $x_1$ is supported on $x_1 \asymp X_1$, and $\{w_{\pm Y_1}: Y_1 \in (0,\infty)\}$ is a $1$-inert family of functions with support on $\pm t_1 \asymp Y_1$. Then the family $\{ X_1^{-1} w_{\pm Y_1}(t_1)  \widehat{w}_T(t_1,x_2,\cdots,x_d) \ : \ (T, \pm Y_1) \in \mathcal{F} \times \pm(0,\infty)\}$ is $X$-inert. Furthermore if $Y_1 \gg p^\varepsilon \frac{X}{X_1}$, then for any $A>0$, we have
\begin{equation}
    \label{eq:inertfourier}
     X_1^{-1} w_{\pm Y_1}(t_1)  \widehat{w}_T(t_1,x_2,\cdots,x_d) \ll_{\varepsilon, A} p^{-A}.
\end{equation}
\end{prop}

We can similarly describe the Mellin transform of such functions as well. We state the following result (\cite[Lemma 4.2]{KY21}).

\begin{prop}
\label{prop:MellinInert}
Suppose that $\{ w_T(x_1,x_2,\cdots,x_d) : T \in \mathcal{F}\}$ is a family of $X$-inert functions such that $x_1$ is supported on $x_1 \asymp X_1$. Let
\begin{equation}
    \label{eq:inertmellin}
     \widetilde{w}_T(s,x_2,\cdots,x_d) = \int_{0}^{\infty} w_T(x,x_2,\cdots,x_d) x^s \frac{dx}{x}.
\end{equation}
Then we have $ \widetilde{w}_T(s,x_2,\cdots,x_d) = {X_1}^s W_T(s,x_2,\cdots,x_d)$ where $W_T(s,\cdots)$ is a family of $X$-inert functions in $x_2,\cdots,x_d$, which is entire in $s$, and has rapid decay for $\vert{\Im(s)} \rvert \gg {X_1}^{1+\varepsilon}$.
\end{prop}

The following is a restatement of Lemma 3.1 in \cite{KPY19}.
\begin{prop}
\label{prop:stationeryphaseinert}
Suppose that $w$ is an $X$-inert function, with compact support on $[Z,2Z]$, so that $w^{(j)}(t) \ll (Z/X)^{-j}$. And suppose that $\phi$ is smooth and satisfies $\phi^{(j)}(t) \ll \frac{Y}{Z^j}$ for some $Y/X^2 \geq R \geq 1$, and for all $t \in [Z,2Z]$. Let
\begin{equation}
    I = \int_{-\infty}^{\infty} w(t)e^{i\phi(t)} dt.
\end{equation}
We then have
\begin{enumerate}
    \item[(a)] If $\vert {\phi'(t)} \rvert \gg \frac{Y}{Z}$ for all $t \in [Z,2Z]$, then $I \ll_A ZR^{-A}$ for $A$ arbitrary large.
    \item[(b)] If $\vert \phi''(t) \rvert \gg \frac{Y}{Z^2}$ for all $t \in [Z,2Z]$, and there exists a (necessarily unique) $t_0 \in \mathbb{R}$ such that $\phi'(t_0) = 0$, then 
\begin{equation}
\label{eq:statphaselemmabound}
    I = \frac{e^{i\phi(t_0)}}{\sqrt{\phi''(t_0)}} F(t_0) + O_A(Z R^{-A}),
\end{equation}
where $F$ is an $X$-inert function supported on $t_0 \asymp Z$.
\end{enumerate}
\end{prop}

The previous result has a natural generalization (Main Theorem in \cite{KPY19}). 

\begin{prop}

Suppose that $w$ is an $X$-inert function in $t_1,\cdots, t_d $, supported on $t_1 \asymp Z$ and $t_i \asymp X_i$ for $i=2,3,\dots,d$. Suppose that the smooth function $\phi$ satisfies
\begin{equation}
  \frac{\partial^{a_1+a_2+\cdots+a_d}}{\partial t_1^{a_1} \partial t_2^{a_2}\dots\partial t_d^{a_d}} \phi(t_1,t_2,\dots,t_d) \ll_C \frac{Y}{Z^{a_1}}\frac{1}{X_2^{a_2}\dots X_d^{a_d}},
\end{equation}
 for all $a_1,a_2,\dots,a_d \in \mathbb{N}$.
\\Now, suppose that $\phi'(t_1,t_2,\dots,t_d) \gg \frac{Y}{Z^2}$ (here $\phi''$ and $\phi'$ denote the derivatives of $\phi$ with respect to $t_1$), for all $t_1,t_2,\dots,t_d$ in the support of $w$, and there exists a (necessarily unique) $t_0 \in \mathbb{R}$ such that $\phi'(t_0)=0$. Suppose also that $\frac{Y}{X^2} \geq R \geq 1$, then
\begin{equation}
    I = \int_{\mathbb{R}} e^{i\phi(t_1,\dots,t_d) w(t_1,\dots,t_d)}dt_1 = \frac{Z}{\sqrt{Y}}e^{i\phi(t_0,t_2,\dots,t_d)}W(t_2,\dots,t_d) + O_A(ZR^{-A}),
\end{equation}
for some $X$-inert function $W$, and $A$ can be taken to be arbitrarily large. The implied constant depends on $A$ and $C_{\mathcal{F}}$.

\end{prop}

\subsection{Smooth weight functions}\label{weight}

If a smooth $p^{\varepsilon}$ ($\varepsilon$ is sufficiently small) weight function $V$ has fixed compact support in  $(0, \infty)^n$, we call it \emph{nice}. Actually, all the smooth weight functions $V$ we consider in our paper are \emph{nice}.

For a nice function  $V$, we may
separate variables in $V(x_1, \ldots, x_n)$ by first inserting a redundant function $V(x_1) \cdots V(x_n)$ that is 1 on the support of $V$ and then  applying  Mellin inversion
\begin{displaymath}
\begin{split}
&V(x_1, \ldots, x_n) = V(x_1, \ldots, x_n)V(x_1) \cdots V(x_n) \\
&= \int_{\Re(s_1)=0} \cdots \int_{\Re(s_n)=0} \widehat{V}(s_1, \ldots, s_n) \Big( V(x_1) \cdots V(x_n) x_1^{-s_1} \cdots x_n^{-s_n}\Big) \frac{ds_1 \cdots ds_n}{(2\pi i)^n}.
\end{split}
\end{displaymath}
Here we can truncate the vertical integrals at height $|\Im (s)| \ll p^{\varepsilon}$ at the cost of a negligible error.

\section{Proof of Theorem \ref{triple}}

\subsection{Setup and reduction of the estimation}

We recall that $f$ and $g$ are holomorphic forms of weight $k$, level 1, $h$ is a holomorphic form of weight $2k$, level $p^3$. By \cite[Corollary 1.1]{Hu17}, the triple product formula gives
\begin{equation}\label{triple product formula}
   L\left(\frac{1}{2}, f \times g \times h\right) \ll_{\varepsilon} p^{3+\varepsilon} \cdot \left|\frac{1}{p^3}\int_{\Gamma_{0}(p^3) \bs \mathbb{H}} f(z)y^{k/2} g(z) \vert_{\left(\begin{smallmatrix} p^3 & \\ & 1 \end{smallmatrix} \right)} y^{k/2} \overline{h(z)} y^k \frac{dxdy}{y^2} \right|^2, 
\end{equation}
by noting that $p^{-\varepsilon} \ll_{\varepsilon} L_{\text{fin}}(1,f,\text{sym}^2), L_{\text{fin}}(1,g,\text{sym}^2), L_{\text{fin}}(1,h,\text{sym}^2) \ll_{\varepsilon} p^{\varepsilon}$, and $f$, $g$, $h$ are $L^2$-normalized with respect to the normalized $L^2$-norm for weight $k$ holomorphic cusp forms
$$  \frac{1}{[\SL_2(\BZ):\Gamma_0(p^3)]}\cdot \int_{\Gamma_{0}(p^3) \bs \mathbb{H}} \left|F(z)\right|^2 y^k \frac{dxdy}{y^2},$$
where $[\SL_2(\BZ):\Gamma_0(p^3)]=p^3\cdot (1+1/p)=p^3+p^2$.

Since $g$ is a holomorphic cusp form, we can write $g$ as a fixed linear combination of Poincar\'e series. We also note that the first $m$ Poincar\'e series with $1 \leq m \leq \dim S_k(1)$ span the space $S_k(1)$ and $k$ is fixed, without loss of generality, we can assume that
$$
g(z) = P_{1,1}(z)= \sum_{\gamma \in \Gamma_{\infty} \bs \SL_2(\mathbb{Z})} (cz+d)^{-k} e(\gamma z),
$$
which is the level 1 Poincar\'e  series with shift 1. 
Since $\vert$ is the slash operator for weight $k$, we have $$ g(z) \vert_{\left( \begin{smallmatrix} p^3 & \\ & 1 \end{smallmatrix} \right)} y^{k/2} =p^{3k / 2} g(p^3 z) y^{k/2}.$$
In order to prepare our forthcoming unfolding step, we slightly substitute the Poincare series $p^{3k / 2} g(p^3 z)$ with $p^{3k / 2} P_{p^3,  p^3}(z)$. Here
\begin{equation} \label{g expansion}
   p^{3k / 2} P_{p^3, p^3}(z):=p^{3k / 2} \sum_{\gamma \in \Gamma_{\infty} \backslash \Gamma_0(p^3)}(c z+d)^{-k} e(p^3 \gamma z)  
\end{equation}
is the level $p^3$ Poincar\'e  series with shift $p^3$. Before unfolding, we compute the $L^2$-norm of the difference of the two Poincare series 
\begin{equation}\label{difference of the two Poincare series}
  \int_{\Gamma_0(p^3) \bs \mathbb{H}} \left \vert p^{3k/2}g(p^3z)- p^{3k/2}P_{p^3,p^3}(z)\right \rvert^2 y^k\frac{dxdy}{y^2}. 
\end{equation} 
By definition, \eqref{difference of the two Poincare series} is equal to
$$ \vert \vert p^{3k/2} g(p^3z) \rvert \rvert_{L^2}^2+ \vert \vert p^{3k/2} P_{p^3,p^3}(z) \rvert \rvert_{L^2}^2-2 \Re \langle p^{3k/2} g(p^3z), p^{3k/2} P_{p^3,p^3}(z)\rangle,$$
where $\langle , \rangle$ is the inner product defined on the $L^2$-space.

By \cite[Corollary 3.4]{Iwan97}, we have
\begin{equation}\label{L_2 norm 1}
\begin{aligned}
\vert \vert p^{3k/2} g(p^3z) \rvert \rvert_{L^2}^2 &= [\SL_2(\BZ): \Gamma_0(p^3)] \cdot \frac{\Gamma(k-1)}{(4\pi)^{k-1}} \cdot \left( 1+ 2\pi i^{-k} \sum_{c \geq 1} \frac{S(1,1,c)}{c} J_{k-1}\left( \frac{4\pi}{c}\right)\right) \\
&= (p^3+p^2) \cdot \frac{\Gamma(k-1)}{(4\pi)^{k-1}} \cdot \left( 1+ 2\pi i^{-k} \sum_{c \geq 1} \frac{S(1,1,c)}{c} J_{k-1}\left( \frac{4\pi}{c}\right)\right).
\end{aligned}
\end{equation}
Similarly, we have
\begin{equation}
\vert \vert p^{3k/2} P_{p^3,p^3}(z) \rvert \rvert_{L^2}^2= p^3 \cdot \frac{\Gamma(k-1)}{(4\pi)^{k-1}} \cdot \left( 1+ 2\pi i^{-k} \sum_{c \geq 1} \frac{S(p^3,p^3,cp^3)}{cp^3} J_{k-1}\left( \frac{4\pi}{c}\right)\right).
\end{equation}
By exclusion-inclusion, we note that $\sum_{c \geq 1} \frac{S(p^3,p^3,cp^3)}{cp^3} J_{k-1}\left( \frac{4\pi}{c}\right)$ equals
\begin{equation}
\begin{aligned}
   &  \sum_{(c,p)=1} \frac{S(p^3,p^3,cp^3)}{cp^3} J_{k-1}\left( \frac{4\pi}{c}\right)+ \sum_{c=c_1p} \frac{S(p^3,p^3,c_1p^4)}{c_1p^4} J_{k-1}\left( \frac{4\pi}{c_1p}\right) \\
  =& \frac{S(0,0,p^3)}{p^3} \cdot \sum_{(c,p)=1} \frac{S(1,1,c)}{c} J_{k-1}\left( \frac{4\pi}{c}\right) + \sum_{c=c_1p} \frac{S(p^3,p^3,c_1p^4)}{c_1p^4} J_{k-1}\left( \frac{4\pi}{c_1p}\right) \\
  =& \frac{S(0,0,p^3)}{p^3} \cdot \sum_{c \geq 1} \frac{S(1,1,c)}{c} J_{k-1}\left( \frac{4\pi}{c}\right) - \frac{S(0,0,p^3)}{p^3} \cdot \sum_{c=c_2p} \frac{S(1,1,c_2p)}{c_2p} J_{k-1}\left( \frac{4\pi}{c_2p}\right) \\
  &+  \sum_{c=c_1p} \frac{S(p^3,p^3,c_1p^4)}{c_1p^4} J_{k-1}\left( \frac{4\pi}{c_1p}\right).
\end{aligned}
\end{equation}
Since $k>50$ is a fixed even number, applying \cite[Lemma 3.6, Lemma 3.7]{Kim24} and taking the prime number $p$ sufficiently large (For example, $p \gg e^k$), we obtain
$$ \frac{S(0,0,p^3)}{p^3} \cdot \sum_{c=c_2p} \frac{S(1,1,c_2p)}{c_2p} J_{k-1}\left( \frac{4\pi}{c_2p}\right) \ll p^{-30}, \; \sum_{c=c_1p} \frac{S(p^3,p^3,c_1p^4)}{c_1p^4} J_{k-1}\left( \frac{4\pi}{c_1p}\right) \ll p^{-30} $$
by noting that $J_{k-1}(1/p) \ll p^{-(k-1)} < p^{-40}$.
Hence, we have
\begin{equation}\label{L_2 norm 2}
\begin{aligned}
& \vert \vert p^{3k/2} P_{p^3,p^3}(z) \rvert \rvert_{L^2}^2 \\ =& p^3 \cdot \frac{\Gamma(k-1)}{(4\pi)^{k-1}} \cdot \left( 1+ 2\pi i^{-k} \cdot \frac{S(0,0,p^3)}{p^3} \cdot \sum_{c \geq 1} \frac{S(1,1,c)}{c} J_{k-1}\left( \frac{4\pi}{c}\right)+O(p^{-30})\right) \\
=& p^3 \cdot \frac{\Gamma(k-1)}{(4\pi)^{k-1}} \cdot \left( 1+ 2\pi i^{-k} \cdot \left( 1-\frac{1}{p} \right) \cdot \sum_{c \geq 1} \frac{S(1,1,c)}{c} J_{k-1}\left( \frac{4\pi}{c}\right)+O(p^{-30})\right).
\end{aligned}
\end{equation}
Moreover, by \cite[Theorem 3.3]{Iwan97}, we have  
\begin{equation}\label{the inner product}
\langle p^{3k/2} g(p^3z), p^{3k/2} P_{p^3,p^3}(z)\rangle= p^3 \cdot \frac{\Gamma(k-1)}{(4\pi)^{k-1}} \cdot \left( 1+ 2\pi i^{-k} \sum_{c \geq 1} \frac{S(1,1,c)}{c} J_{k-1}\left( \frac{4\pi}{c}\right)\right).
\end{equation}
By \cite[Corollary 3.4]{Iwan97}, we have
$$ \vert \vert g(z) \rvert \rvert_{L^2}^2= \frac{\Gamma(k-1)}{(4\pi)^{k-1}} \cdot \left( 1+ 2\pi i^{-k} \sum_{c \geq 1} \frac{S(1,1,c)}{c} J_{k-1}\left( \frac{4\pi}{c}\right)\right).$$
Therefore, $ \langle p^{3k/2} g(p^3z), p^{3k/2} P_{p^3,p^3}(z)\rangle \in \BR$.

Hence, by \eqref{L_2 norm 1}, \eqref{L_2 norm 2} and \eqref{the inner product}, we have the upper bound
\begin{equation}\label{difference of the two Poincare series upper bound1}
 \int_{\Gamma_0(p^3) \bs \mathbb{H}} \left \vert p^{3k/2}g(p^3z)- p^{3k/2}P_{p^3,p^3}(z)\right \rvert^2 y^k\frac{dxdy}{y^2} \ll p^2.   
\end{equation}
Applying Cauchy--Schwarz inequality we have the upper bound by \cite{Xia07}:
\begin{equation}\label{difference of the two Poincare series upper bound2}
\begin{aligned}
 & \left \vert \int_{\Gamma_0(p^3) \bs \mathbb{H}} \vert f(z) \rvert y^{k/2}\left \vert p^{3k/2}g(p^3z)- p^{3k/2}P_{p^3,p^3}(z)\right \rvert y^{k/2} \vert h(z) \rvert y^k \frac{dxdy}{y^2} \right \rvert^2 \\
 \ll & \left \vert \sup_{z \in \mathbb{H}} y^{k/2} \vert f(z) \rvert \right \rvert^2 \times  \int_{\Gamma_0(p^3) \bs \mathbb{H}} \left \vert p^{3k/2}g(p^3z)- p^{3k/2}P_{p^3,p^3}(z)\right \rvert^2 y^k\frac{dxdy}{y^2} \\
 & \times \int_{\Gamma_0(p^3) \bs \mathbb{H}} \vert h(z) \rvert^2y^{2k} \frac{dxdy}{y^2} \\
 \ll_{\varepsilon}  &  k^{1/2+\varepsilon} \cdot p^2 \cdot p^{3+\varepsilon} \\
 \ll_{\varepsilon} &  p^{5+\varepsilon}.
 \end{aligned}
\end{equation}

By \eqref{triple product formula},  \eqref{difference of the two Poincare series upper bound2} and the unfolding technique, we obtain 
\begin{equation}
\begin{aligned}
& L(1 / 2, f \times g \times h) \\
 \ll_{\varepsilon} & p^{3+\varepsilon} \left|\frac{1}{p^3} \int_{\Gamma_{0}(p^3) \backslash \mathbb{H}} \sum_{\gamma \in \Gamma_{\infty} \backslash \Gamma_0(p^3)}p^{3k / 2}(c z+d)^{-k}  e(p^3 \gamma z) f(z) \overline{h(z)} y^{2k}  d \mu(z)\right|^2 + p^{2+\varepsilon}\\
 =& p^{3+\varepsilon} \left|\frac{p^{3k / 2}}{p^3} \int_{\Gamma_{0}(p^3) \backslash \mathbb{H}} \sum_{\gamma \in \Gamma_{\infty} \backslash \Gamma_0(p^3)}|c z+d|^{-4k}  e(p^3 \gamma z) f(\gamma z) \overline{h(\gamma z)}   y^{2k}\frac{d x d y}{y^2}\right|^2 + p^{2+\varepsilon}\\
=& p^{3+\varepsilon} \left|\frac{p^{3k / 2}}{p^3} \int_{\Gamma_{0}(p^3) \backslash \mathbb{H}} \sum_{\gamma \in \Gamma_{\infty} \backslash \Gamma_0(p^3)}(\Im(\gamma z))^{2k}  e(p^3 \gamma z) f(\gamma z) \overline{h(\gamma z)}   \frac{d x d y}{y^2}\right|^2 + p^{2+\varepsilon}\\
=& p^{3+\varepsilon} \left|\frac{p^{3k / 2}}{p^3} \int_{\Gamma_{\infty} \backslash \mathbb{H}} f(z) \overline{h(z)} e(p^3 x) e^{-2 \pi p^3 y} y^{2k}d \mu(z)\right|^2 + p^{2+\varepsilon}.
\end{aligned}
\end{equation}
Inserting the Fourier expansion, we obtain
\begin{equation}\label{triple product expansion}
\begin{aligned}
L(1 / 2, f \times g \times h) & \ll_{\varepsilon} p^{3+\varepsilon} \left|\frac{p^{3k / 2}}{p^3} \sum_{\substack{m=n+p^3 \\
n \geq 1}} \lambda_{f}(n) \lambda_{h}(m) \frac{n^{(k-1) / 2 } \cdot m^{k-1/2}}{(n+m+p^3)^{2k-1  }}\right|^2 + p^{2+\varepsilon}  \\
& \ll_{\varepsilon} p^{3+\varepsilon} \left|\frac{p^{3k / 2}}{p^3} \sum_{\substack{m=n+p^3 \\
n \geq 1}} \lambda_{f}(n) \lambda_{h}(m) \frac{n^{\frac{k-1}{2} }}{m^{k-1/2  }}\right|^2 + p^{2+\varepsilon}.
\end{aligned}
\end{equation}

By dyadic partition of unity and opening the square, it suffices to prove a nontrivial upper bound for
\begin{equation} \label{double}
 \frac{1}{p^3} \cdot \sum_{n_1, n_2} \lambda_f(n_1)\lambda_f(n_2)\lambda_h(n_1+p^3)\lambda_h(n_2+p^3)V \left( \frac{n_1}{N} \right)V \left( \frac{n_2}{N} \right),
\end{equation}
where $V$ is a smooth weight ($1$-inert) function for $n_1$ and $n_2$ with support on  $[1,2]$,
and $p^{3-\varepsilon} \leq N \leq p^{3+\varepsilon}$ since we pick weight $k>k_{\varepsilon}$ which is fixed and sufficiently large. A trivial estimation of the above double sums on $n_1,n_2$ will cover the convexity bound $p^{3+\varepsilon}$.

\begin{rmk}
Because of the extra term $p^{2+\varepsilon}$, we cannot prove a short first moment estimation for the triple product $L$-function in the cubic level aspect, which can be used to deduce the strong subconvexity bound for the triple product $L$-function. However, since $p^{2+\varepsilon}$ is small enough, we can still apply the refined Petersson trace formula in \cite{PWZ20} and deduce the strong subconvexity bound because of the nonnegativity of the double sum in \eqref{double}.
\end{rmk}

\subsection{Applying Petersson trace formula }
We recall $p^{-\varepsilon} \ll_{\varepsilon} L_{\text{fin}}(1,h,\text{sym}^2) \ll_{\varepsilon} p^{\varepsilon}$ and the nonnegativity for the double summation in \eqref{double}. In order to prove our main Theorem \ref{triple}, it suffices to show that

\begin{equation}
S:=\frac{1}{p^3}\sum_{n_1 \geq 1}\sum_{n_2 \geq 1}
\lambda_{f}(n_1) \lambda_{f}(n_2) V \left( \frac{n_1}{N} \right)V \left( \frac{n_2}{N} \right) \sum_{h \in \CA(2k,M,p^3)} 
 \frac{\lambda_{h}(n_1+p^3)\lambda_{h}(n_2+p^3)}{L_{\text{fin}}(1,h,\text{sym}^2)} \ll_{\varepsilon} p^{2+2\theta+\varepsilon},
\end{equation}
where $1 \leq M \leq p-1$ is a fixed and unique positive integer. 
Since $h$ is a newform of level $p^3$, we note that $\lambda_h(n_1+p^3)=\lambda_h(n_2+p^3)=0$ if $p \vert n_1$ and $p \vert n_2$.

Applying the refined Petersson trace formula in Theorem \ref{peter1} (see also \cite[Proposition 3.1]{PWZ20}), we obtain 
$$
S = \mathcal{D} + \frac{(-1)^k(2k-1)}{\pi} \mathcal{O},
$$
where 
$$
\mathcal{D} := \frac{1}{p^3}\sum_{\substack{n_1 \geq 1 \\ (n_1,p)=1}}\sum_{\substack{n_2 \geq 1 \\ (n_2,p)=1}}
\lambda_{f}(n_1) \lambda_{f}(n_2) V \left( \frac{n_1}{N} \right)V \left( \frac{n_2}{N} \right)\delta(n_1+p^3,n_2+p^3) \frac{(2k-1)p^2}{2\pi^2} \ll_{\varepsilon} p^{2+\varepsilon},
$$
and 
$$
\begin{aligned}
\mathcal{O} := \frac{1}{ p^3} & \sum_{\substack{n_1 \geq 1 \\ (n_1,p)=1}}\sum_{\substack{n_2 \geq 1 \\ (n_2,p)=1}}
\lambda_{f}(n_1) \lambda_{f}(n_2)  V \left( \frac{n_1}{N} \right)V \left( \frac{n_2}{N} \right) \sum_{c \geq 1} \frac{A_{p,M}(c)}{c} \\
&\cdot S(n_1+p^3,n_2+p^3, p^2c) \cdot J_{2k-1} \left( \frac{4\pi \sqrt{(n_1+p^3)(n_2+p^3)}}{p^2c} \right).
\end{aligned}
$$
We note that the summation in $\CO$ is in terms of $(n_1 n_2, p)=1$. Since $k>k_{\varepsilon}$ and is fixed, we may truncate the $c$-sum to $c \ll p^{1+\varepsilon}$. Moreover, we can insert a smooth (nice) weight function $V(c/C)$ for some $1 \leq C \leq p^{1+\varepsilon}$. It suffices to prove that
$$
\begin{aligned}
\frac{1}{ p^3} & \sum_{\substack{n_1 \geq 1 \\ (n_1,p)=1}}\sum_{\substack{n_2 \geq 1 \\ (n_2,p)=1}}
\lambda_{f}(n_1) \lambda_{f}(n_2)  \sum_{c \geq 1} \frac{A_{p,M}(c)}{c} \cdot V \left (\frac{n_1}{N} \right)V \left (\frac{n_2}{N} \right)V \left (\frac{c}{C} \right)\\
&\cdot S(n_1+p^3,n_2+p^3, p^2c) \cdot J_{2k-1} \left( \frac{4\pi \sqrt{(n_1+p^3)(n_2+p^3)}}{p^2c} \right) \\
& \ll_{\varepsilon} p^{2+2\theta+\varepsilon}.
\end{aligned}
$$

By exclusion-inclusion on the $c$-sums, it suffices to rewrite $\mathcal{O}$ and bound the following two terms:
$$
\begin{aligned}
\mathcal{O}_1:= \frac{1}{ p^3} & \sum_{\substack{n_1 \geq 1 \\ (n_1,p)=1}}\sum_{\substack{n_2 \geq 1 \\ (n_2,p)=1}}
\lambda_{f}(n_1) \lambda_{f}(n_2)  \sum_{\substack{c \geq 1 \\ (c,p)=1}} \frac{1}{c} \cdot V \left (\frac{n_1}{N} \right)V \left (\frac{n_2}{N} \right)V \left (\frac{c}{C} \right)\\
&\cdot S(n_1+p^3,n_2+p^3, p^2c) \cdot J_{2k-1} \left( \frac{4\pi \sqrt{(n_1+p^3)(n_2+p^3)}}{p^2c} \right),
\end{aligned}
$$
and
$$
\begin{aligned}
\mathcal{O}_2:= \frac{1}{ p^3} & \sum_{\substack{n_1 \geq 1 \\ (n_1,p)=1}}\sum_{\substack{n_2 \geq 1 \\ (n_2,p)=1}}
\lambda_{f}(n_1) \lambda_{f}(n_2)  \sum_{\substack{c \geq 1 \\ p \vert c}} \frac{1}{c} \cdot V \left (\frac{n_1}{N} \right)V \left (\frac{n_2}{N} \right)V \left (\frac{c}{C} \right)\\
&\cdot S(n_1+p^3,n_2+p^3, p^2c) \cdot J_{2k-1} \left( \frac{4\pi \sqrt{(n_1+p^3)(n_2+p^3)}}{p^2c} \right).
\end{aligned}
$$
By \cite[Proposition 5.1]{michel0}, we have $\mathcal{O}_2 \ll_{\varepsilon} p^{2+\varepsilon}$. 

By exclusion-inclusion on the $n_1$ and $n_2$-sums, it suffices to rewrite $\mathcal{O}_1$ and bound the following four terms:
$$
\begin{aligned}
\mathcal{O}_3:= \frac{1}{ p^3} & \sum_{\substack{n_1 \geq 1 }}\sum_{\substack{n_2 \geq 1 }}
\lambda_{f}(n_1) \lambda_{f}(n_2)  \sum_{\substack{c \geq 1 }} \frac{1}{c} \cdot V \left (\frac{n_1}{N} \right)V \left (\frac{n_2}{N} \right)V \left (\frac{c}{C} \right)\\
&\cdot S(n_1+p^3,n_2+p^3, p^2c) \cdot J_{2k-1} \left( \frac{4\pi \sqrt{(n_1+p^3)(n_2+p^3)}}{p^2c} \right),
\end{aligned}
$$
$$
\begin{aligned}
\mathcal{O}_4:= \frac{1}{ p^3} & \sum_{\substack{n_1 \geq 1 \\ (n_1,p)=1}}\sum_{\substack{n_2 \geq 1 \\ p \vert n_2}}
\lambda_{f}(n_1) \lambda_{f}(n_2)  \sum_{\substack{c \geq 1 }} \frac{1}{c} \cdot V \left (\frac{n_1}{N} \right)V \left (\frac{n_2}{N} \right)V \left (\frac{c}{C} \right)\\
&\cdot S(n_1+p^3,n_2+p^3, p^2c) \cdot J_{2k-1} \left( \frac{4\pi \sqrt{(n_1+p^3)(n_2+p^3)}}{p^2c} \right),
\end{aligned}
$$
$$
\begin{aligned}
\mathcal{O}_5:= \frac{1}{ p^3} & \sum_{\substack{n_1 \geq 1 \\ p \vert n_1}}\sum_{\substack{n_2 \geq 1 \\ (n_2,p)=1}}
\lambda_{f}(n_1) \lambda_{f}(n_2)  \sum_{\substack{c \geq 1 }} \frac{1}{c} \cdot V \left (\frac{n_1}{N} \right)V \left (\frac{n_2}{N} \right)V \left (\frac{c}{C} \right)\\
&\cdot S(n_1+p^3,n_2+p^3, p^2c) \cdot J_{2k-1} \left( \frac{4\pi \sqrt{(n_1+p^3)(n_2+p^3)}}{p^2c} \right),
\end{aligned}
$$
and
$$
\begin{aligned}
\mathcal{O}_6:= \frac{1}{ p^3} & \sum_{\substack{n_1 \geq 1 \\ p \vert n_1}}\sum_{\substack{n_2 \geq 1 \\ p \vert n_2}}
\lambda_{f}(n_1) \lambda_{f}(n_2)  \sum_{\substack{c \geq 1 }} \frac{1}{c} \cdot V \left (\frac{n_1}{N} \right)V \left (\frac{n_2}{N} \right)V \left (\frac{c}{C} \right)\\
&\cdot S(n_1+p^3,n_2+p^3, p^2c) \cdot J_{2k-1} \left( \frac{4\pi \sqrt{(n_1+p^3)(n_2+p^3)}}{p^2c} \right).
\end{aligned}
$$
If $(n_1,p)=1$ and $p \vert n_2$, we have $S(n_1+p^3,n_2+p^3,cp^2)=0$ since $p^2 \vert cp^2$, $(n_1+p^3,p)=1$, $p \vert n_2+p^3$. With the same reason, we have $S(n_1+p^3,n_2+p^3,cp^2)=0$ if $p \vert n_1$ and $(n_2,p)=1$. Hence, we have $\mathcal{O}_4=\mathcal{O}_5=0$.

\begin{lem}
We have $$ \mathcal{O}_6 \ll_{\varepsilon} p^{2+\varepsilon}.$$
\end{lem}

\begin{proof}
We can give the estimation by applying the Kuznetsov trace formula and the spectral large sieve inequality. We consider two cases: (1) $1 \leq C \leq p^{1-\varepsilon}$ and (2) $p^{1-\varepsilon} \leq C \leq p^{1+\varepsilon}$.

If $1 \leq C \leq p^{1-\varepsilon}$, then we have $\sqrt{(n_1+p^3)(n_2+p^3)}/{(p^2c)} \gg p^{\varepsilon}$. For $t \gg p^{\varepsilon}$, the Bessel function has an oscillatory behavior and satisfies
$$ J_{2k-1}(t)=\frac{1}{\sqrt{t}} \cdot (e^{it}W_+(t)+e^{-it}W_{-}(t)),$$
where $W_+$ and $W_{-}$ satisfy the same derivative bounds as $1$-inert functions. We insert a $1$-inert function $W \left( 4 \pi \sqrt{(n_1+p^3)(n_2+p^3)}/{(p^2c)} \right)$, which is supported on $1 \leq X<x <2X$ with $p/C \leq X \leq p^{1+\varepsilon}/C$. Up to a $p^{\varepsilon}$ factor, we can remove the weight function $V \left( c/C \right)$ by a Mellin inversion (see \cite[Lemma 2.10]{Khan22}, Section \ref{subsec:inertfunctions} and \ref{weight}). Hence, it suffices to prove that
$$
\begin{aligned}
\widetilde{\mathcal{O}}& :=\frac{1}{ p^3}  \sum_{\substack{n_1 \geq 1 \\ p \vert n_1}}\sum_{\substack{n_2 \geq 1 \\ p \vert n_2}}
\lambda_{f}(n_1) \lambda_{f}(n_2) \left(\frac{p}{C} \right)^{-1/2}  V \left (\frac{n_1}{N} \right)V \left (\frac{n_2}{N} \right) \sum_{\substack{c \geq 1 }} \frac{1}{c}  \\
&\cdot S(n_1+p^3,n_2+p^3, p^2c) \cdot e \left( \frac{\pm 2 \sqrt{(n_1+p^3)(n_2+p^3)}}{p^2c} \right) \cdot W \left( \frac{4\pi \sqrt{(n_1+p^3)(n_2+p^3)}}{p^2c} \right) \\
& \ll_{\varepsilon} p^{2+\varepsilon}.
\end{aligned}
$$
For the inner $c$-sum, we can apply Kuznetsov trace formula (Theorem \ref{kuznetsov}). This will transfer the $c$-sum as follows:
$$ \sum_{\substack{g \in B(p^2,1) \\ t_g \ll T}} \frac{1}{\vert \vert g \rvert \rvert^2} \frac{1}{\cosh({\pi t_g})} \cdot \widetilde{\phi}(t_g)\lambda_g(n_1+p^3)\overline{\lambda_g(n_2+p^3)}+(\cdots), $$
where $(\cdots)$ denotes the continuous spectrum contribution and the holomorphic contribution, which are similar to deal with. Applying Lemma \ref{transf-bounds} (c), we have $\widetilde{\phi}(t_g) \ll (p/C)^{-1/2+\varepsilon}$ and $\widetilde{\phi}(t_g) \ll p^{-A}$ for arbitray large $A>0$ if $t_g \gg (p/C)^{1/2+\varepsilon}$. Hence, we can always restrict the spectral parameters in the Kuznetsov trace formula (Theorem \ref{kuznetsov}) and the spectral large sieve inequality (Theorem \ref{largesievethm}) to $k, \vert t_g \rvert , \vert t \rvert \ll (p/C)^{1/2+\varepsilon}$. 

Applying Cauchy--Schwarz inequality, Lemma \ref{transf-bounds} (c) and Theorem \ref{largesievethm} (Here we have $T=(p^{1+\varepsilon}/C)^{1/2}$ and $\sum_{p \vert n_1} \vert \lambda_f(n_1) \rvert^2 V \left( n_1/N\right) \ll_{\varepsilon} p^{2+\varepsilon}$), it suffices to bound $1/p^3 \cdot T^{-2} \cdot \sqrt{S_1}\cdot \sqrt{S_2} \cdot p^2$, where
$$ S_1:= \sum_{\substack{g \in B(p^2,1) \\ t_g \ll T}} \frac{1}{\vert \vert g \rvert \rvert^2} \frac{1}{\cosh({\pi t_g})} \cdot \left \vert \sum_{\substack{n_1 \geq p^3+1 \\ p \vert n_1}} V\left(\frac{n_1-p^3}{N}\right)\lambda_f(n_1-p^3)\lambda_g(n_1)\right \rvert^2+(\cdots), $$
and
$$ S_2:= \sum_{\substack{g \in B(p^2,1) \\ t_g \ll T}} \frac{1}{\vert \vert g \rvert \rvert^2} \frac{1}{\cosh({\pi t_g})} \cdot \left \vert \sum_{\substack{n_2 \geq p^3+1 \\ p \vert n_2}} V\left(\frac{n_2-p^3}{N}\right)\lambda_f(n_2-p^3)\lambda_g(n_2)\right \rvert^2+(\cdots),$$
where $(\cdots)$ denotes the continuous spectrum contribution and the holomorphic contributions, which are similar to deal with. 
We have the following upper bound for $\widetilde{\mathcal{O}}$:
$$ \widetilde{\mathcal{O}} \ll_{\varepsilon} \frac{1}{p^3} \cdot \left( \frac{p^{1+\varepsilon}}{C}\right)^{-1} \cdot \left( \frac{p^{1+\varepsilon}}{C}+ \frac{p^{3+\varepsilon}}{p^2}\right) \cdot p^{2+\varepsilon} \cdot p^2 \ll_{\varepsilon} p^{2+\varepsilon}.$$ 

If $p^{1-\varepsilon} \leq C \leq p^{1+\varepsilon}$, then we have $p^{-\varepsilon} \ll \sqrt{(n_1+p^3)(n_2+p^3)}/{(p^2c)} \ll p^{\varepsilon}$. In this region, the Bessel function is not oscillatory. For $p^{-\varepsilon} \ll t \ll p^{\varepsilon}$, We can write $J_{2k-1}(t)=t^{2k-1}W(t)$, where $t^jW^{(j)}(t) \ll T_0$ with $T_0 \ll p^{\varepsilon}$. This is the same derivative bound satisfied by a $T_0$-inert function. We insert a $1$-inert function $W \left( 4\pi \sqrt{(n_1+p^3)(n_2+p^3)}/{(p^2c)} \right)$, which is supported on $0 < X<x <2X$ with $p^{-\varepsilon} \leq X \leq p^{\varepsilon}$. Up to a $p^{\varepsilon}$ factor, we can remove the weight function $V \left( c/C \right)$ by a Mellin inversion (see Section \ref{subsec:inertfunctions}). Hence, it suffices to prove that
$$
\begin{aligned}
\widetilde{\mathcal{O}}& :=\frac{1}{ p^3}  \sum_{\substack{n_1 \geq 1 \\ p \vert n_1}}\sum_{\substack{n_2 \geq 1 \\ p \vert n_2}}
\lambda_{f}(n_1) \lambda_{f}(n_2)   V \left (\frac{n_1}{N} \right)V \left (\frac{n_2}{N} \right) \sum_{\substack{c \geq 1 }} \frac{1}{c}  \\
&\cdot S(n_1+p^3,n_2+p^3, p^2c)  \cdot W \left( \frac{4\pi \sqrt{(n_1+p^3)(n_2+p^3)}}{p^2c} \right) \\
& \ll_{\varepsilon} p^{2+\varepsilon}.
\end{aligned}
$$
For the inner $c$-sums, we can apply Theorem \ref{kuznetsov}. Applying Cauchy--Schwarz, Lemma \ref{transf-bounds} (a) and Theorem \ref{largesievethm} (Here we have $T=p^{\varepsilon}$ and $\sum_{p \vert n_1} \vert \lambda_f(n_1) \rvert^2 V \left( n_1/N\right) \ll_{\varepsilon} p^{2+\varepsilon}$), we have the following upper bound for $\widetilde{\mathcal{O}}$:
$$ \widetilde{\mathcal{O}} \ll_{\varepsilon} \frac{1}{p^3} \cdot  \left( p^{\varepsilon}+ \frac{p^{3+\varepsilon}}{p^2}\right) \cdot p^{2+\varepsilon} \cdot p^2 \ll_{\varepsilon} p^{2+\varepsilon}.$$
Finally, for the possible exceptional eigenvalues, from Lemma \ref{transf-bounds} (b), we note that $\tilde{\phi}(t_g) \ll_{\varepsilon} p^{\varepsilon}$ since $p^{-\varepsilon} \ll X/Z \ll p^{\varepsilon}$. Now, we apply Theorem \ref{kuznetsov}, Cauchy--Schwarz inequality and Theorem \ref{largesievethm} (Here we have $\vert t_g \rvert \ll 1$ and $\sum_{p \vert n_1} \vert \lambda_f(n_1) \rvert^2 V \left( n_1/N\right) \ll_{\varepsilon} p^{2+\varepsilon}$), we have the following upper bound for $\widetilde{\mathcal{O}}$:
$$ \widetilde{\mathcal{O}} \ll_{\varepsilon} \frac{1}{p^3} \cdot  \left( p^{\varepsilon}+ \frac{p^{3+\varepsilon}}{p^2}\right) \cdot p^{2+\varepsilon} \cdot p^2 \ll_{\varepsilon} p^{2+\varepsilon}.$$
\end{proof}

Next, it suffices to prove that $\mathcal{O}_3 \ll_{\varepsilon} p^{2+2\theta+\varepsilon}$.

\subsection{Voronoi summation and stationary phase}

Opening the Kloosterman sum and exchanging the order of summations, we get
$$
\begin{aligned}
\mathcal{O}_3 = &\frac{1}{ p^3} \sum_{n_1 }\lambda_{f}(n_1) V \left (\frac{n_1}{N} \right)
\sum_{\substack{c \geq 1 }} V \left (\frac{c}{C} \right) \frac{1}{c} 
\sideset{}{^\star}\sum_{a\bmod{p^2c}}
 e\left(\frac{an_1 }{p^2 c}\right)e\left(\frac{(a+\overline{a})p^3}{cp^2}\right)\\
&\cdot \sum_{n_2 } \lambda_{f}(n_2) e\left(\frac{\overline{a} n_2}{p^2 c}\right) 
V\left( \frac{n_2}{N} \right) J_{2k-1}\left( \frac{4 \pi \sqrt{(n_1+p^3)(n_2+p^3)}}{p^2 c} \right).
\end{aligned}
$$

$$
\begin{aligned}
& \sum_{n_1 }\lambda_{f}(n_1) V \left (\frac{n_1}{N} \right)
\sum_{\substack{c \geq 1 }} V \left (\frac{c}{C} \right) \frac{1}{c} 
\sideset{}{^\star}\sum_{a\bmod{p^2c}}
 e\left(\frac{an_1 }{p^2 c}\right)e\left(\frac{(a+\overline{a})p^3}{cp^2}\right)\\
&\times  \sum_{n_2 } \lambda_{f}(n_2) e\left(\frac{\overline{a} n_2}{p^2 c}\right) 
V\left( \frac{n_2}{N} \right) J_{2k-1}\left( \frac{4 \pi \sqrt{(n_1+p^3)(n_2+p^3)}}{p^2 c} \right) \\
\ll_{\varepsilon} & p^{5+2\theta+\varepsilon}. 
\end{aligned}
$$

Applying the Voronoi summation formula to the above $n_2$-sum, we obtain
$$
\frac{1}{p^2c} \sum_{n_2} \lambda_{f}(n_2) e\left(\frac{-a n_2}{p^2c}\right) \Phi_{1}\left(\frac{n_2}{p^4c^2}\right),
$$
where 
$$
\begin{aligned}
    \Phi_{1}\left(\frac{n_2}{p^4c^2}\right)= 2\pi i^k \int_{0}^{\infty} & V\left(\frac{y}{N} \right) J_{2k-1}\left( \frac{4 \pi \sqrt{(n_1+p^3)(y+p^3)}}{p^2 c} \right) J_{k-1}\left(\frac{4\pi\sqrt{n_2y}}{p^2c}\right) dy. 
\end{aligned}
$$
By making a change of variable, we have
$$
\Phi_{1}\left(\frac{n_2}{p^4c^2}\right)= 2\pi i^k N \int_{0}^{\infty}  V\left(y \right) J_{2 k-1}\left( \frac{4 \pi \sqrt{(n_1+p^3)(yN+p^3)}}{p^2 c} \right) J_{k-1}\left(\frac{4\pi\sqrt{n_2 N y}}{p^2 c}\right) d y. 
$$

We consider two different cases: (1) $1 \leq C  \ll p^{1-\varepsilon}$ and (2) $p^{1-\varepsilon} \ll C \ll p^{1+\varepsilon}$.

If $p^{1-\varepsilon} \ll C \ll p^{1+\varepsilon}$, we may assume that $N \gg p^3$. The case for $N \ll p^3$ is similar. It is known that $p^{-\varepsilon} \ll \sqrt{(n_1+p^3)(yN+p^3)}/(p^2c) \ll p^{\varepsilon}$. In this region, the Bessel function is not oscillatory. For $p^{-\varepsilon} \ll t \ll p^{\varepsilon}$, We can write $J_{2k-1}(t)=t^{2k-1}W(t)$, where $t^jW^{(j)}(t) \ll T_0$ with $T_0 \ll p^{\varepsilon}$. This is the same derivative bound satisfied by a $T_0$-inert function. If $n_2 \ll p^{3-\varepsilon}$, since the weight $k$ is a sufficiently large fixed even number, the contribution is negligible. If $n_2 \gg p^{3+\varepsilon}$, it is known that $\sqrt{n_2N}/(p^2c) \gg p^{\varepsilon}$. By integration by parts, the contribution from this part is negligible. We can now restrict the $n_2$ sum to $p^{3-\varepsilon} \ll n_2 \ll p^{3+\varepsilon}$. By a dyadic partition for the new $n_2$-sum, in order to understand $\Phi_1\left( n_2/(p^4c^2)\right)$, it suffices to consider
$$ N \cdot V \left( \frac{n_2}{N_2}\right),$$
where $V$ is a $p^{\varepsilon}$-inert function supported on $n_2 \asymp N_2$, and $p^{3-\varepsilon} \ll N_2 \ll p^{3+\varepsilon}$. We note that $p^{-\varepsilon} \ll N/(p^2C) \ll p^{\varepsilon}$. Hence, for $p^{1-\varepsilon} \ll C \ll p^{1+\varepsilon}$, it suffices to prove that
\begin{equation}\label{nonos}
\begin{aligned}
& \mathcal{K}_1:=\sum_{n_1 \geq 1} \sum_{n_2 \geq 1}\lambda_{f}(n_1) \lambda_f(n_2) V \left (\frac{n_1}{N} \right) V \left (\frac{n_2}{N_2} \right) 
 \sum_{\substack{c \geq 1 }} V \left (\frac{c}{C} \right) \frac{1}{c} \cdot S(p^3,n_1-n_2+p^3,p^2c)  \\ 
\ll_{\varepsilon} & p^{5+2\theta+\varepsilon}. 
\end{aligned}
\end{equation}

If $1 \leq C \ll p^{1-\varepsilon}$, then we have $\sqrt{(n_1+p^3)(yN+p^3)}/(p^2c) \gg p^{\varepsilon}$. For $t \gg p^{\varepsilon}$, the Bessel function has an oscillatory behavior and satisfies
$$ J_{k-1}(t)=\frac{1}{\sqrt{t}} \cdot (e^{it}W_+(t)+e^{-it}W_{-}(t)),$$
where $W_+$ and $W_{-}$ satisfy the same derivative bounds as $1$-inert functions. We can also assume that $\sqrt{n_2N}/(p^2c) \gg p^{\varepsilon}$. Otherwise, integration by part gives the negligible integral. Hence, it suffices to consider the following integral:
$$
N \cdot \int_{0}^{\infty}  V\left(y \right) J_{2 k-1}\left( \frac{4 \pi \sqrt{(n_1+p^3)(yN+p^3)}}{p^2 c} \right) J_{k-1}\left(\frac{4\pi\sqrt{n_2 N y}}{p^2 c}\right) d y. 
$$
We may assume that $N \gg p^3$. The case for $N \ll p^3$ is similar. Applying $ J_{k-1}(t)=\frac{1}{\sqrt{t}} \cdot (e^{it}W_+(t)+e^{-it}W_{-}(t))$, it suffices to consider the following integral:
$$
\begin{aligned}
 N \cdot & \int_{0}^{\infty}  \widetilde{V}\left(y \right) 
 e\left( \frac{\pm 2 \sqrt{(n_1+p^3)(yN+p^3)}}{p^2 c} \right) \cdot 
 e\left(\frac{\pm 2\sqrt{n_2 N y}}{p^2 c}\right) d y, 
\end{aligned}
$$
where
$$ 
\begin{aligned}
\widetilde{V}\left(y \right)&:= V(y)W_{\pm}\left( \frac{ \sqrt{(n_1+p^3)(yN+p^3)}}{p^2 c} \right)W_{\pm} \left( \frac{\sqrt{n_2 N y}}{p^2 c}\right) \\
& \times \left( \frac{\sqrt{(n_1+p^3)(yN+p^3)}}{p^2c}\right)^{-1/2}\left( \frac{\sqrt{n_2Ny}}{p^2c}\right)^{-1/2}.
\end{aligned}
$$
Hence, $\widetilde{V}$ is a $p^{\varepsilon}$-inert function in terms of $y$ and is compactly support on the interval $[1,2]$. Moreover, we only need to consider the case that the sign inside the exponential function is opposite. Otherwise, by integration by parts, the above integral is negligible.
It suffices to consider the following integral:
\begin{equation}\label{y-integral}
  N \cdot \int_{0}^{\infty}  \widetilde{V}\left(y \right) 
e\left( \frac{ 2 \sqrt{(n_1+p^3)(yN+p^3)}}{p^2 c} -\frac{ 2\sqrt{n_2 N y}}{p^2 c} \right)  d y. 
\end{equation}
Let $$ f(y):= \frac{ 2 \sqrt{(n_1+p^3)(yN+p^3)}}{p^2 c} -\frac{ 2\sqrt{n_2 N y}}{p^2 c},$$
then we have
$$ f'(y)= \frac{ N (n_1+p^3)^{1/2}(yN+p^3)^{-1/2}}{p^2 c} -\frac{ \sqrt{n_2 N }y^{-1/2}}{p^2 c}.$$
It is known that $n_2 \asymp N$. Otherwise, we have the lower bound $f'(y) \gg N/(p^2C) \gg p^{\varepsilon}$, in which case the integral is negligible due to Proposition \ref{prop:stationeryphaseinert} (a). We set $f'(y)=0$ and get the unique stationary point $$ t_0= \frac{n_2 p^3}{(n_1-n_2+p^3)N}.$$
By direct computation, we can assume that $n_1-n_2+p^3 \asymp p^3$ and $t_0 \in [4/5,11/5]$. Otherwise, we have the lower bound $f'(y) \gg p/C \gg p^{\varepsilon}$ for all $y \in [1/2]$, and we know that the above integral is negligible by Proposition \ref{prop:stationeryphaseinert} (a).

We have 
$$ f(t_0)=\frac{2 \sqrt{(n_1-n_2+p^3)p^3}}{p^2c} \asymp \frac{p}{C}$$
and
$$ f''(t_0) \asymp \frac{p}{C}.$$
Applying stationary phase method (Proposition \ref{prop:stationeryphaseinert} (b)) to  \eqref{y-integral}, we obtain 
$$ p^2C \cdot \frac{e(f(t_0))}{\sqrt{f''(t_0)}} \cdot G(t_0)+ O_A\left( \left( \frac{p}{C} \right)^{-A} \right),$$
where $G$ is a $p^{\varepsilon}$-inert function compactly supported on $(1/10,10)$. 
Hence, it is reduced to estimate
\begin{equation}\label{os}
\begin{aligned}
& \mathcal{K}_2:=\sum_{n_1 \geq 1} \sum_{n_2 \geq 1}\lambda_{f}(n_1) \lambda_f(n_2) V \left (\frac{n_1}{N} \right) V \left (\frac{n_2}{N} \right) 
  \sum_{\substack{c \geq 1 }} V \left (\frac{c}{C} \right) \left( \frac{p}{C}\right)^{-1/2} \\ 
  & \times e\left( \frac{2 \sqrt{(n_1-n_2+p^3)p^3}}{p^2c}\right) \cdot G \left( \frac{n_2p^3}{(n_1-n_2+p^3)N} \right) \cdot \frac{1}{c} \cdot S(p^3,n_1-n_2+p^3,p^2c)  \\ 
\ll_{\varepsilon} & p^{5+2\theta+\varepsilon}. 
\end{aligned}
\end{equation}

Next, we will estimate $\mathcal{K}_1$ and $\mathcal{K}_2$.
We write $r:=n_1-n_2$. The case $r=0$ will be considered later. We further assume that $r>0$. The case $r<0$ is similar to the case $r>0$. Hence, for $p^{1-\varepsilon} \ll C \ll p^{1+\varepsilon}$, \eqref{nonos} yields
\begin{equation}   \label{nonos1}
\begin{aligned}
& \sum_{n_1 \geq 1} \sum_{r \geq 1}\lambda_{f}(n_1) \lambda_f(n_1-r) V \left (\frac{n_1}{N} \right) V \left (\frac{n_1-r}{N_2} \right) 
 \sum_{\substack{c \geq 1 }} V \left (\frac{c}{C} \right) \frac{1}{c} \cdot S(p^3,r+p^3,p^2c).   
\end{aligned}
\end{equation}
For a fixed $r \geq 1$, we define
$$ F(r):=\sum_{c \geq 1} V\left( \frac{c}{C} \right) \frac{1}{c} \cdot S\left( p^3,r+p^3,p^2c\right).$$
We may also trivially truncate the $r$-sum in  \eqref{nonos1} by $1 \leq r \ll N$. Up to a $p^{\varepsilon}$ factor, we may restrict the $r$-sum to $N_3 \leq r < 2N_3$, where $1 \leq N_3 \ll N \ll p^{3+\varepsilon}$. If $r<0$, we may consider $\widetilde{r}:=-r>0$. Similarly, we may restrict the $\widetilde{r}$-sum to $\widetilde{N_3} \leq \widetilde{r}< 2 \widetilde{N_3}$ and $\widetilde{N_4} \leq p^3-\widetilde{r} < 2\widetilde{N_4}$ or $\widetilde{N_4} \leq \widetilde{r}-p^3 < 2\widetilde{N_4}$.

Since $1 \leq r \ll p^{3+\varepsilon}$, we can split \eqref{nonos1} into four different cases: (1) $(r,p)=1$; (2) $p \vert \vert r$; (3) $p^2 \vert \vert r$ and (4) $p^3 \vert \vert r$.
The four cases are the following:
\begin{equation*}   
\begin{aligned}
E_0:= \sum_{n_1 \geq 1} \sum_{(r,p)=1}\lambda_{f}(n_1) \lambda_f(n_1-r) V \left (\frac{n_1}{N} \right) V \left (\frac{n_1-r}{N_2} \right) 
 \sum_{\substack{c \geq 1 }} V \left (\frac{c}{C} \right) \frac{1}{c} \cdot S(p^3,r+p^3,p^2c);  \\ 
\end{aligned}
\end{equation*}
\begin{equation} \label{E_i}
\begin{aligned}
E_i:= \sum_{n_1 \geq 1} \sum_{p^i \vert \vert r}\lambda_{f}(n_1) \lambda_f(n_1-r) V \left (\frac{n_1}{N} \right) V \left (\frac{n_1-r}{N_2} \right) 
 \sum_{\substack{c \geq 1 }} V \left (\frac{c}{C} \right) \frac{1}{c} \cdot S(p^3,r+p^3,p^2c),  \\ 
\end{aligned}
\end{equation}
for $i=1, 2, 3$.

For the special case $r=0$, we will reduce to the sum
\begin{equation}   
\begin{aligned}
E:=& \sum_{n_1 \geq 1} \lambda_{f}(n_1) \lambda_f(n_1) V \left (\frac{n_1}{N} \right) V \left (\frac{n_1}{N_2} \right) 
 \sum_{\substack{c \geq 1 }} V \left (\frac{c}{C} \right) \frac{1}{c} \cdot S(p^3,p^3,p^2c).  \\ 
\end{aligned}
\end{equation}
It suffices to prove that $E_0,E_1,E_2,E_3,E \ll_{\varepsilon} p^{5+2\theta+\varepsilon}$.

For $1 \leq C \ll p^{1-\varepsilon}$,  \eqref{os} yields
\begin{equation}  \label{os1}
\begin{aligned}
& \sum_{n_1 \geq 1} \sum_{r \geq 1}\lambda_{f}(n_1) \lambda_f(n_1-r) V \left (\frac{n_1}{N} \right) V \left (\frac{n_1-r}{N} \right) 
  \sum_{\substack{c \geq 1 }} V \left (\frac{c}{C} \right) \left( \frac{p}{C}\right)^{-1/2} \\ 
  & \times e\left( \frac{2 \sqrt{(r+p^3)p^3}}{p^2c}\right) \cdot G \left( \frac{(n_1-r)p^3}{(r+p^3)N} \right) \cdot \frac{1}{c} \cdot S(p^3,r+p^3,p^2c).  \\ 
\end{aligned}
\end{equation}
For \eqref{os1} and a fixed $r \geq 1$, we define
$$ \widetilde{F}(r):= \sum_{c \geq 1} V \left( \frac{c}{C}\right) \left( \frac{p}{C} \right)^{-1/2} e\left( \frac{2 \sqrt{(r+p^3)p^3}}{p^2c}\right) \cdot \frac{1}{c} \cdot S(p^3, r+p^3,p^2c).$$
We may also trivially truncate the $r$-sum in  \eqref{os1} by $1 \leq r \ll \min{(N,p^3)}$. Up to a $p^{\varepsilon}$ factor, we may restrict the $r$-sum to $N_4 \leq r < 2N_4$, where $1 \leq N_4 \ll N \ll p^{3+\varepsilon}$. If $r<0$, we may consider $\widetilde{r}:=-r>0$. Similarly, we may restrict the $\widetilde{r}$-sum to $\widetilde{N_5} \leq \widetilde{r}< 2 \widetilde{N_5}$ and $\widetilde{N_6} \leq p^3-\widetilde{r} < 2\widetilde{N_6}$. 


We can remove the $p^{\varepsilon}$-inert weight function $G$ by a Mellin inversion (see \cite[Lemma 2.10]{Khan22}, Section \ref{subsec:inertfunctions} and \ref{weight}). Hence, up to a $p^{\varepsilon}$ factor, it suffices to consider the summation
\begin{equation}  \label{os2}
\begin{aligned}
& \sum_{n_1 \geq 1} \sum_{r \geq 1}\lambda_{f}(n_1) \lambda_f(n_1-r) V \left (\frac{n_1}{N} \right) V \left (\frac{n_1-r}{N} \right) \\
  & \times \sum_{\substack{c \geq 1 }} V \left (\frac{c}{C} \right) \left( \frac{p}{C}\right)^{-1/2}  
  \cdot  e\left( \frac{2 \sqrt{(r+p^3)p^3}}{p^2c}\right) \cdot  \frac{1}{c} \cdot S(p^3,r+p^3,p^2c),  \\ 
\end{aligned}
\end{equation}
where $V$ is certain $p^{\varepsilon}$-inert function.

Since $1 \leq r \ll p^{3+\varepsilon}$, we can split the summation \eqref{os2} into four different cases: (1) $(r,p)=1$; (2) $p \vert \vert r$; (3) $p^2 \vert \vert r$ and (4) $p^3 \vert \vert r$.
The four cases are the following:
\begin{equation*}  
\begin{aligned}
F_0:=& \sum_{n_1 \geq 1} \sum_{(r,p)=1}\lambda_{f}(n_1) \lambda_f(n_1-r) V \left (\frac{n_1}{N} \right) V \left (\frac{n_1-r}{N} \right) \\
  & \times \sum_{\substack{c \geq 1 }} V \left (\frac{c}{C} \right) \left( \frac{p}{C}\right)^{-1/2}  
  \cdot  e\left( \frac{2 \sqrt{(r+p^3)p^3}}{p^2c}\right) \cdot  \frac{1}{c} \cdot S(p^3,r+p^3,p^2c);  \\ 
\end{aligned}
\end{equation*}
\begin{equation}  \label{F_i}
\begin{aligned}
F_i:=& \sum_{n_1 \geq 1} \sum_{p^{i} \vert \vert r}\lambda_{f}(n_1) \lambda_f(n_1-r) V \left (\frac{n_1}{N} \right) V \left (\frac{n_1-r}{N} \right) \\
  & \times \sum_{\substack{c \geq 1 }} V \left (\frac{c}{C} \right) \left( \frac{p}{C}\right)^{-1/2}  
  \cdot  e\left( \frac{2 \sqrt{(r+p^3)p^3}}{p^2c}\right) \cdot  \frac{1}{c} \cdot S(p^3,r+p^3,p^2c), 
  \end{aligned}
\end{equation}
for $i=1,2,3$.

For the special case $r=0$, we will reduce to the sum
\begin{equation}  
\begin{aligned}
F:=& \sum_{n_1 \geq 1} \lambda_{f}(n_1) \lambda_f(n_1) V \left (\frac{n_1}{N} \right) V \left (\frac{n_1}{N} \right) \\
  & \times \sum_{\substack{c \geq 1 }} V \left (\frac{c}{C} \right) \left( \frac{p}{C}\right)^{-1/2}  
  \cdot  e\left(  \frac{ 2p}{c}\right) \cdot  \frac{1}{c} \cdot S(p^3,p^3,p^2c).  \\ 
\end{aligned}
\end{equation}
It suffices to prove that $F_0,F_1,F_2,F_3,F \ll_{\varepsilon} p^{5+2\theta+\varepsilon}$.

\subsection{$L^2$-norm for $F(r)$ and $\widetilde{F}(r)$} \label{l2}

As the preparation for the final application of the spectral large sieve inequality (Section \ref{finalkuz}), we have to compute the $L^2$-norm for $F(r)$ and $\widetilde{F}(r)$. The main ingredient for the estimation is the spectral theory and application for the Kuznetsov trace formula (Theorem \ref{kuznetsov}) backwards.

For a fixed $r \geq 1$, we recall
$$ F(r):=\sum_{c \geq 1} V\left( \frac{c}{C} \right) \frac{1}{c} \cdot S\left( p^3,r+p^3,p^2c\right),$$
for $p^{1-\varepsilon} \ll C \ll p^{1+\varepsilon}$. For $1 \leq C \ll p^{1-\varepsilon}$, we recall
$$ \widetilde{F}(r):= \sum_{c \geq 1} V \left( \frac{c}{C}\right) \left( \frac{p}{C} \right)^{-1/2} e\left( \frac{2 \sqrt{(r+p^3)p^3}}{p^2c}\right) \cdot \frac{1}{c} \cdot S(p^3, r+p^3,p^2c).$$

In this subsection, we will give an upper bound for the following $L^2$-norm by the method of Kuznetsov trace formula and the spectral large sieve inequality. The target terms are the following:
$$ \sum_{(r,p)=1} \vert F(r) \rvert^2,\;\;\; \sum_{(r,p)=1} \vert \widetilde{F}(r) \rvert^2,\;\;\;\sum_{p \vert \vert r} \vert F(r) \rvert^2,\;\;\; \sum_{p \vert \vert r} \vert \widetilde{F}(r) \rvert^2,$$
$$ \sum_{p^2 \vert \vert r} \vert F(r) \rvert^2,\;\;\; \sum_{p^2 \vert \vert r} \vert \widetilde{F}(r) \rvert^2,\;\;\;\sum_{p^3 \vert \vert r} \vert F(r) \rvert^2,\;\;\; \sum_{p^3 \vert \vert r} \vert \widetilde{F}(r) \rvert^2.$$
We note that the $r$-sum can be restricted to $1 \leq r \ll p^{3+\varepsilon}$.
If $(r,p)=1$, since $p \vert p^3, p \nmid r+p^3,p^2 \vert p^2c$, by the property of Kloosterman sums, we have $S(p^3,r+p^3,p^2c)=0$. Hence, $E_0=F_0 \equiv 0$ and $\sum_{(r,p)=1} \vert F(r) \rvert^2= \sum_{(r,p)=1} \vert \widetilde{F}(r) \rvert^2 \equiv 0$. If $p \vert c$ (this is equivalent to $p \vert \vert c$ since $c \ll p^{1+\varepsilon}$), we also recall that $S(p^3,r+p^3,p^2c)=0$ if $p \vert \vert r$, which can be deduced by the factorization of Kloosterman sums and the exact formula for the Ramanujan sums.

Now we compute $\sum_{p \vert \vert r} \vert \widetilde{F}(r) \rvert^2$, and the treatment for $L^2$-norm $\sum_{p \vert \vert r} \vert F(r) \rvert^2$ is similar. Since $p \vert \vert r$, we write $r:=pr_1$ and $(r_1,p)=1$. By definition, we have
\begin{equation}
\begin{aligned}
 \widetilde{F}(r)& = \sum_{(c,p)=1} V \left( \frac{c}{C}\right) \left( \frac{p}{C} \right)^{-1/2} e\left( \frac{2 \sqrt{(r+p^3)p^3}}{p^2c}\right) \cdot \frac{1}{c} \cdot S(p^3, r+p^3,p^2c) \\
 &+ \sum_{c=pc_1,\,c_1 \geq 1} V \left( \frac{pc_1}{C}\right) \left( \frac{p}{C} \right)^{-1/2} e\left( \frac{2 \sqrt{(r+p^3)p^3}}{p^3c_1}\right) \cdot \frac{1}{pc_1} \cdot S(p^3, r+p^3,p^3c_1).
 \end{aligned}
 \end{equation}
Since $c \ll p^{1+\varepsilon}$, we have $c_1 \ll p^{\varepsilon}$ and $(c_1,p)=1$. By the factorization of Kloosterman sums and the exact formula for the Ramanujan sums, we have $S(p^3,r+p^3,p^3c_1) \equiv 0$ since $\mu(p^2)=\mu(p^3)=0$ by the definition of the Mobius function. This gives that $$ \widetilde{F}(r) = \sum_{(c,p)=1} V \left( \frac{c}{C}\right) \left( \frac{p}{C} \right)^{-1/2} e\left( \frac{2 \sqrt{(r+p^3)p^3}}{p^2c}\right) \cdot \frac{1}{c} \cdot S(p^3, r+p^3,p^2c). $$
Since $(c,p)=1$, we have $ S(p^3,r+p^3,p^2c)=-p\cdot S(1,r_1+p^2,c).$ Hence, $\widetilde{F}(r)=-p \cdot \widetilde{H}(r_1)$, where
$$ \widetilde{H}(r_1) := \sum_{(c,p)=1} V \left( \frac{c}{C}\right) \left( \frac{p}{C} \right)^{-1/2} e\left( \frac{2 \sqrt{r_1+p^2}}{c}\right) \cdot \frac{1}{c} \cdot S(1, r_1+p^2,c). $$
We rewrite $\widetilde{H}(r_1)$ as follows:
\begin{equation}
\begin{aligned}
\widetilde{H}(r_1) &= \sum_{c \geq 1} V \left( \frac{c}{C}\right) \left( \frac{p}{C} \right)^{-1/2} e\left( \frac{2 \sqrt{r_1+p^2}}{c}\right) \cdot \frac{1}{c} \cdot S(1, r_1+p^2,c) \\
& - \sum_{c=pc_1,\,c_1 \geq 1} V \left( \frac{pc_1}{C}\right) \left( \frac{p}{C} \right)^{-1/2} e\left( \frac{2 \sqrt{r_1+p^2}}{pc_1}\right) \cdot \frac{1}{pc_1} \cdot S(1, r_1+p^2,pc_1).
\end{aligned}
\end{equation}
By the Weil's bound for Kloosterman sums and $c_1 \ll p^{\varepsilon}$, we deduce that the square of the second term in the above equation has the upper bound $p^{-1+\varepsilon}$. 
For the first term in the above equation, we insert a $1$-inert function $W\left( 4 \pi \sqrt{r_1+p^2}/c\right)$, which is supported on $1 \leq X <x <2X$ with $X \asymp p/C$. Up to a $p^{\varepsilon}$ factor, we can remove the weight function $V\left(c/C \right)$ by a Mellin inversion. Now we apply the Kuznetsov trace formula (Theorem \ref{kuznetsov}) to
$$ \left( \frac{p}{C} \right)^{-1/2} \cdot \sum_{c \geq 1}   e\left( \frac{2 \sqrt{r_1+p^2}}{c}\right) \cdot W \left( \frac{4\pi \sqrt{r_1+p^2}}{c} \right) \cdot \frac{1}{c} \cdot S(1, r_1+p^2,c),$$
which equals $(p/C)^{-1/2} \cdot S_1$, with
$$ S_1:= \sum_{\substack{g \in B(1,1) \\ t_g \ll T}} \frac{1}{\vert \vert g \rvert \rvert^2} \frac{1}{\cosh({\pi t_g})} \cdot \widetilde{\phi}(t_g)\lambda_g(1)\overline{\lambda_g(r_1+p^2)}+(\cdots), $$
where $(\cdots)$ denotes the continuous spectrum contribution and the holomorphic contribution, which are similar to deal with. Here we have $T=(p^{1+\varepsilon}/C)^{1/2}$. We note that $\widetilde{\phi}(t_g) \ll (p/C)^{-1/2+\varepsilon}$ and $\lambda_g(r_1+p^2) \ll p^{2\theta+\varepsilon}$.
Applying the spectral large sieve inequality (Theorem \ref{largesievethm}) (Weyl law), Cauchy--Schwarz inequality, Lemma \ref{transf-bounds} (c) 
and the Rankin--Selberg bound, we have an upper bound for $\widetilde{H}(r_1)$, which is given by $T^{-2}\cdot T^2\cdot p^{2\theta+\varepsilon} \ll p^{2\theta+\varepsilon}$. Therefore,
$$ \sum_{(r_1,p)=1} \vert \widetilde{H}(r_1) \rvert^2 \ll \frac{p^{3+\varepsilon}}{p} \cdot T^{-4+\varepsilon} \cdot T^2 \cdot T^2 \cdot p^{4\theta+\varepsilon}+ \frac{p^{3+\varepsilon}}{p} \cdot p^{-1+\varepsilon} \ll p^{2+4\theta+\varepsilon}.$$

Similarly, we also have $\sum_{(r_1,p)=1} \vert H(r_1) \rvert^2 \ll p^{2+4\theta+\varepsilon}$.

For the case $p^2 \vert \vert r$, we write $r:=p^2r_2$ and $(r_2,p)=1$. By definition, we have
\begin{equation}
\begin{aligned}
 \widetilde{F}(r)& = \sum_{(c,p)=1} V \left( \frac{c}{C}\right) \left( \frac{p}{C} \right)^{-1/2} e\left( \frac{2 \sqrt{(r+p^3)p^3}}{p^2c}\right) \cdot \frac{1}{c} \cdot S(p^3, r+p^3,p^2c) \\
 &+ \sum_{c=pc_1,\,c_1 \geq 1} V \left( \frac{pc_1}{C}\right) \left( \frac{p}{C} \right)^{-1/2} e\left( \frac{2 \sqrt{(r+p^3)p^3}}{p^3c_1}\right) \cdot \frac{1}{pc_1} \cdot S(p^3, r+p^3,p^3c_1).
 \end{aligned}
 \end{equation}
 Since $c \ll p^{1+\varepsilon}$, we have $c_1 \ll p^{\varepsilon}$ and $(c_1,p)=1$. If $(c,p)=1$, we have $S(p^3,r+p^3,p^2c)=(p^2-p)\cdot S(p,r_2+p,c).$ If $(c_1,p)=1$, we have $S(p^3,r+p^3,p^3c_1)=-p^2 \cdot S(1,\overline{p}r_2+1,c_1)$. We write $\widetilde{F}(r)=(p^2-p)\cdot \widetilde{H}(r_2)$, where
 \begin{equation}
 \begin{aligned}
  \widetilde{H}(r_2) := & \sum_{(c,p)=1} V \left( \frac{c}{C}\right) \left( \frac{p}{C} \right)^{-1/2} e\left( \frac{2 \sqrt{p(r_2+p)}}{c}\right) \cdot \frac{1}{c} \cdot S(p, r_2+p,c) \\
   -& \frac{p^2}{p^2-p} \cdot \sum_{c=pc_1,\,c_1 \geq 1} V \left( \frac{pc_1}{C}\right) \left( \frac{p}{C} \right)^{-1/2} e\left( \frac{2 \sqrt{p(r_2+p)}}{pc_1}\right) \cdot \frac{1}{pc_1} \cdot S(1, \overline{p}r_2+1,c_1).
  \end{aligned}
  \end{equation}
We rewrite $\widetilde{H}(r_2)$ as follows:
\begin{equation}
\begin{aligned}
\widetilde{H}(r_2) &= \sum_{c \geq 1} V \left( \frac{c}{C}\right) \left( \frac{p}{C} \right)^{-1/2} e\left( \frac{2 \sqrt{p(r_2+p)}}{c}\right) \cdot \frac{1}{c} \cdot S(p, r_2+p,c) \\
-& \sum_{c=pc_1,\,c_1 \geq 1} V \left( \frac{pc_1}{C}\right) \left( \frac{p}{C} \right)^{-1/2} e\left( \frac{2 \sqrt{p(r_2+p)}}{pc_1}\right) \cdot \frac{1}{pc_1} \cdot S(p, r_2+p,pc_1) \\
 -& \frac{p^2}{p^2-p} \cdot \sum_{c=pc_1,\,c_1 \geq 1} V \left( \frac{pc_1}{C}\right) \left( \frac{p}{C} \right)^{-1/2} e\left( \frac{2 \sqrt{p(r_2+p)}}{pc_1}\right) \cdot \frac{1}{pc_1} \cdot S(1, \overline{p}r_2+1,c_1).
\end{aligned}
\end{equation}
By the Weil's bound for Kloosterman sums and $c_1 \ll p^{\varepsilon}$, we deduce that the square of the second and the third term in the above equation have the upper bound $p^{-1+\varepsilon}$ and $p^{-2+\varepsilon}$. For the first term in the above equation, we insert a $1$-inert function $W\left( 4 \pi \sqrt{p(r_2+p)}/c\right)$, which is supported on $1 \leq X <x <2X$ with $X \asymp p/C$. Up to a $p^{\varepsilon}$ factor, we can remove the weight function $V\left(c/C \right)$ by a Mellin inversion. Now we apply the Kuznetsov trace formula (Theorem \ref{kuznetsov}) to
$$ \left( \frac{p}{C} \right)^{-1/2} \cdot \sum_{c \geq 1}   e\left( \frac{2 \sqrt{p(r_2+p)}}{c}\right) \cdot W \left( \frac{4\pi \sqrt{p(r_2+p)}}{c} \right) \cdot \frac{1}{c} \cdot S(p, r_2+p,c). $$
Applying the spectral large sieve inequality (Theorem \ref{largesievethm}) (Weyl law), Cauchy--Schwarz inequality, Lemma \ref{transf-bounds} (c) (Here we have $T=(p^{1+\varepsilon}/C)^{1/2}$), and the Rankin--Selberg bound, we have
$$ \sum_{(r_2,p)=1} \vert \widetilde{H}(r_2) \rvert^2 \ll \frac{p^{3+\varepsilon}}{p^2} \cdot T^{-4+\varepsilon} \cdot T^2 \cdot T^2 \cdot p^{4\theta+\varepsilon}+ \frac{p^{3+\varepsilon}}{p^2} \cdot (p^{-1+\varepsilon}+p^{-2+\varepsilon}) \ll p^{1+4\theta+\varepsilon}.$$

Similarly, we also have $\sum_{(r_2,p)=1} \vert H(r_2) \rvert^2 \ll p^{1+4\theta+\varepsilon}$.

For the case $p^3  \vert r$, we write $r:=p^3r_3$ and $(r_3,p)=1$. By definition, we have
\begin{equation}
\begin{aligned}
 \widetilde{F}(r)& = \sum_{(c,p)=1} V \left( \frac{c}{C}\right) \left( \frac{p}{C} \right)^{-1/2} e\left( \frac{2 \sqrt{(r+p^3)p^3}}{p^2c}\right) \cdot \frac{1}{c} \cdot S(p^3, r+p^3,p^2c) \\
 &+ \sum_{c=pc_1,\,c_1 \geq 1} V \left( \frac{pc_1}{C}\right) \left( \frac{p}{C} \right)^{-1/2} e\left( \frac{2 \sqrt{(r+p^3)p^3}}{p^3c_1}\right) \cdot \frac{1}{pc_1} \cdot S(p^3, r+p^3,p^3c_1).
 \end{aligned}
 \end{equation}
 Since $c \ll p^{1+\varepsilon}$, we have $c_1 \ll p^{\varepsilon}$ and $(c_1,p)=1$. If $(c,p)=1$, we have $S(p^3,r+p^3,p^2c)=(p^2-p)\cdot S(p,p(r_3+1),c).$ If $(c_1,p)=1$, we have $S(p^3,r+p^3,p^3c_1)=(p^3-p^2) \cdot S(1,r_3+1,c_1)$. We write $\widetilde{F}(r)=(p^2-p)\cdot \widetilde{H}(r_3)$, where
 \begin{equation}
 \begin{aligned}
  \widetilde{H}(r_3) := & \sum_{(c,p)=1} V \left( \frac{c}{C}\right) \left( \frac{p}{C} \right)^{-1/2} e\left( \frac{2 \sqrt{p(pr_3+p)}}{c}\right) \cdot \frac{1}{c} \cdot S(p, pr_3+p,c) \\
   +&  \sum_{c=pc_1,\,c_1 \geq 1} V \left( \frac{pc_1}{C}\right) \left( \frac{p}{C} \right)^{-1/2} e\left( \frac{2 \sqrt{p(pr_3+p)}}{pc_1}\right) \cdot \frac{1}{c_1} \cdot S(1, r_3+1,c_1).
  \end{aligned}
  \end{equation}
We rewrite $\widetilde{H}(r_3)$ as follows:
\begin{equation}
\begin{aligned}
\widetilde{H}(r_3) &= \sum_{c \geq 1} V \left( \frac{c}{C}\right) \left( \frac{p}{C} \right)^{-1/2} e\left( \frac{2 \sqrt{p(pr_3+p)}}{c}\right) \cdot \frac{1}{c} \cdot S(p, pr_3+p,c) \\
-& \sum_{c=pc_1,\,c_1 \geq 1} V \left( \frac{pc_1}{C}\right) \left( \frac{p}{C} \right)^{-1/2} e\left( \frac{2 \sqrt{p(pr_3+p)}}{pc_1}\right) \cdot \frac{1}{pc_1} \cdot S(p, pr_3+p,pc_1) \\
 +&  \sum_{c=pc_1,\,c_1 \geq 1} V \left( \frac{pc_1}{C}\right) \left( \frac{p}{C} \right)^{-1/2} e\left( \frac{2 \sqrt{p(pr_3+p)}}{pc_1}\right) \cdot \frac{1}{c_1} \cdot S(1, r_3+1,c_1).
\end{aligned}
\end{equation}
By the Weil's bound for Kloosterman sums and $c_1 \ll p^{\varepsilon}$, we deduce that the square of the second and the third term in the above equation have the upper bound $p^{\varepsilon}$ and $p^{\varepsilon}$. For the first term in the above equation, we insert a $1$-inert function $W\left( 4 \pi \sqrt{p(pr_3+p)}/c\right)$, which is supported on $1 \leq X <x <2X$ with $p^{1-\varepsilon}/C \leq X \leq p^{1+\varepsilon}/C$. Up to a $p^{\varepsilon}$ factor, we can remove the weight function $V\left(c/C \right)$ by a Mellin inversion. Now we apply the Kuznetsov trace formula (Theorem \ref{kuznetsov}) to
$$ \left( \frac{p}{C} \right)^{-1/2} \cdot \sum_{c \geq 1}   e\left( \frac{2 \sqrt{p(pr_3+p)}}{c}\right) \cdot W \left( \frac{4\pi \sqrt{p(pr_3+p)}}{c} \right) \cdot \frac{1}{c} \cdot S(p, pr_3+p,c). $$
Applying the spectral large sieve inequality (Theorem \ref{largesievethm}) (Weyl law), Cauchy--Schwarz inequality, Lemma \ref{transf-bounds} (c) (Here we have $T=(p^{1+\varepsilon}/C)^{1/2}$), and the Rankin--Selberg bound, we have
$$ \sum_{(r_3,p)=1} \vert \widetilde{H}(r_3) \rvert^2 \ll \frac{p^{3+\varepsilon}}{p^3} \cdot T^{-4+\varepsilon} \cdot T^2 \cdot T^2 \cdot p^{4\theta+\varepsilon}+ \frac{p^{3+\varepsilon}}{p^3} \cdot (p^{\varepsilon}+p^{\varepsilon}) \ll p^{4\theta+\varepsilon}.$$

Similarly, we also have $\sum_{(r_3,p)=1} \vert H(r_3) \rvert^2 \ll p^{4\theta+\varepsilon}$. We also note that for $r_3=-1$, the Kloosterman sums reduce to the Ramanujan sum. In this special case, without applying the Kuznetsov trace formula, by the exact formula for the Ramanujan sums (If $(c,p)=1$, we have $S(p^3,0,cp^2)/c=(p^2-p)\mu(c)/c$. If $c=c_1p$ and $(c_1,p)=1$, we have $S(p^3,0,cp^2)/c=(p^2-p)\mu(c_1)/c_1$), we can deduce that $H(-1) \ll_{\varepsilon} p^{\varepsilon}$.

This completes the computation for the $L^2$-norms.

\subsection{The degenerate case: $r=0$}\label{degenerate case r=0}  \label{dege}

For the special case $r=0$, it suffices to prove that $F \ll_{\varepsilon} p^{5+2\theta+\varepsilon}$. Let
$$ \widetilde{F}(0) := \sum_{(c,p)=1} V \left( \frac{c}{C}\right) \left( \frac{p}{C} \right)^{-1/2} e\left( \frac{2p}{c}\right) \cdot \frac{1}{c} \cdot S(p^3, p^3,p^2c). $$
Since $F=\sum_{n_1 \geq 1} \lambda_f(n_1)\lambda_f(n_1)V \left( n_1/N \right) V \left( n_1/N \right) \cdot \widetilde{F}(0)$, by the Rankin--Selberg bound, it suffices to prove that $\widetilde{F}(0) \ll_{\varepsilon} p^{2+2\theta+\varepsilon}$. 

By a similar method as the above discussion for $p^3 \vert r$, it suffices to show that $\widetilde{H}(0) \ll_{\varepsilon} p^{2\theta+\varepsilon}$ (We note that $\widetilde{F}(0)=(p^2-p)\widetilde{H}(0)$), where
\begin{equation}
\begin{aligned}
\widetilde{H}(0) &:= \sum_{c \geq 1} V \left( \frac{c}{C}\right) \left( \frac{p}{C} \right)^{-1/2} e\left( \frac{2 p}{c}\right) \cdot \frac{1}{c} \cdot S(p,p,c) \\
-& \sum_{c=pc_1,\,c_1 \geq 1} V \left( \frac{pc_1}{C}\right) \left( \frac{p}{C} \right)^{-1/2} e\left( \frac{2p}{pc_1}\right) \cdot \frac{1}{pc_1} \cdot S(p, p,pc_1) \\
 +&  \sum_{c=pc_1,\,c_1 \geq 1} V \left( \frac{pc_1}{C}\right) \left( \frac{p}{C} \right)^{-1/2} e\left( \frac{2p}{pc_1}\right) \cdot \frac{1}{c_1} \cdot S(1,1,c_1).
\end{aligned}
\end{equation}
By the Weil's bound for Kloosterman sums and $c_1 \ll p^{\varepsilon}$, we deduce that the square of the second and the third term in the above equation have the upper bound $p^{\varepsilon}$ and $p^{\varepsilon}$. For the first term in the above equation, we insert a $1$-inert function $W\left( 4 \pi p/c\right)$, which is supported on $1 \leq X <x <2X$ with $X \asymp p/C$. Up to a $p^{\varepsilon}$ factor, we can remove the weight function $V\left(c/C \right)$ by a Mellin inversion. Now we apply the Kuznetsov trace formula (Theorem \ref{kuznetsov}) to
$$ \left( \frac{p}{C} \right)^{-1/2} \cdot \sum_{c \geq 1}   e\left( \frac{2 p}{c}\right) \cdot W \left( \frac{4\pi p}{c} \right) \cdot \frac{1}{c} \cdot S(p, p,c),$$
which equals $(p/C)^{-1/2} \cdot S_1$, with
$$ S_1:= \sum_{\substack{g \in B(1,1) \\ t_g \ll T}} \frac{1}{\vert \vert g \rvert \rvert^2} \frac{1}{\cosh({\pi t_g})} \cdot \widetilde{\phi}(t_g)\lambda_g(p)\overline{\lambda_g(p)}+(\cdots), $$
where $(\cdots)$ denotes the continuous spectrum contribution and the holomorphic contribution, which are similar to deal with.
Here we have $T=(p^{1+\varepsilon}/C)^{1/2}$.
Applying the spectral large sieve inequality (Theorem \ref{largesievethm}) (Weyl law), Cauchy--Schwarz inequality and Lemma \ref{transf-bounds} (c), we have
$$ \widetilde{H}(0) \ll_{\varepsilon} T^{-2+\varepsilon} \cdot T^2 \cdot p^{2\theta+\varepsilon}+p^{\varepsilon} \ll_{\varepsilon} p^{2\theta+\varepsilon}$$
by noting that $\widetilde{\phi}(t_g) \ll (p/C)^{-1/2+\varepsilon}$ and $\vert \lambda_g(p) \rvert \leq 2p^{\theta}$. This gives that $\widetilde{F}(0) \ll_{\varepsilon} p^{2+2\theta+\varepsilon}$ and $F \ll_{\varepsilon} p^{5+2\theta+\varepsilon}$.

Similarly, we will also have $F(0) \ll_{\varepsilon} p^{2+2\theta+\varepsilon}$ and $E \ll_{\varepsilon} p^{5+2\theta+\varepsilon}$.

\subsection{Application of the circle method}\label{apply Jutila circle method}

 We now return to the summation $E_i, F_i$ for $i=1,2,3$ in \eqref{E_i} and \eqref{F_i} (Let $n_1=n$, $m=n-r$.) and treat the $n$-sum as a shifted convolution problem:
\begin{equation}\label{apply circle method N'}
   \sum_{n, m} V\Big( \frac{n}{N}, \frac{m}{N_2}\Big) \lambda_{f}(n) \lambda_{f}(m)  \int_0^1 e\big((n-m - r )\alpha\big) d\alpha, 
\end{equation}
or
\begin{equation}\label{apply circle method N}
\sum_{n, m} V\Big( \frac{n}{N}, \frac{m}{N}\Big) \lambda_{f}(n) \lambda_{f}(m)  \int_0^1 e\big((n-m - r )\alpha\big) d\alpha.
\end{equation}
Since $p^{3-\varepsilon} \ll N, N_2 \ll p^{3+\varepsilon}$, we may only treat \eqref{apply circle method N} and the treatment for \eqref{apply circle method N'} is similar.

We choose a gigantic parameter $Q = p^{30}$ and replace the characteristic function on $[0, 1]$  with $I(\alpha)$, defined in Lemma \ref{lem3}. By Cauchy--Schwarz inequality and trivial bounds, this introduces an error at most $p^{6+\varepsilon} \cdot Q^{\eps-1/2}$ which can be neglected. In this way, the $\alpha$-integral becomes
$$ \frac{Q}{2\Lambda} \sum_q V\Big(\frac{q}{Q}\Big)  \sum_{\substack{d\, (\text{{\rm mod }} q)\\ (d, q) = 1}} \int_{-1/Q}^{1/Q} e\left( \Big( \frac{d}{q} + \alpha\Big) (n-m - r ) \right) d\alpha.$$
Since $Q \gg p^{3+\varepsilon}>\max\{m,n,r\}$, it suffices to bound
\begin{equation}\label{beforeVor}
\frac{1}{\Lambda}  \sum_q V\Big(\frac{q}{Q}\Big)\sum_{\substack{d\, (\text{{\rm mod }} q)\\ (d, q) = 1}} \sum_{n, m} V\Big( \frac{n}{N}, \frac{m}{N}\Big) \lambda_{f}(n) \lambda_{f}(m)  e\Big(   \frac{d}{q}  (n-m - r ) \Big),
\end{equation}
where $\Lambda = \sum_q V(q/Q) \phi(q) = Q^{2 + o(1)}$. We recall that this represents the inner $n$-sum in $F_i$ for $i=1,2,3$.

\subsection{Voronoi summation to both variables} \label{twice}

Next, we apply the Voronoi formula to both sums in \eqref{beforeVor}.
For the $m$-sum
$$
\sum_{m} V\Big( \frac{m}{N}\Big)  \lambda_{f}(m)  e\Big( \frac{-d m}{q} \Big),
$$
by Lemma \ref{voronoi} we obtain
$$
\frac{N}{q} \sum_{m}  \lambda_{f}(m) e\Big( \frac{\overline{d} m}{q} \Big) \Phi\left(\frac{mN}{q^2}\right),
$$
where  
$$
\Phi\left(\frac{mN}{q^2}\right)= 2 \pi i^k \int_0^{\infty} V(x) J_{k-1} \left(4 \pi \frac{\sqrt{xmN}}{q }\right) dx .
$$
For $\frac{mN}{q^2} \gg p^{\varepsilon}$, we have 
$$
\Phi\left(\frac{mN}{q^2}\right)= 2 \pi i^k \int_0^{\infty} V(x) \left(\frac{\sqrt{xmN}}{q }\right)^{-1/2} e \left(\frac{\pm 2\sqrt{xmN}}{q }\right)W\left(\frac{\sqrt{xmN}}{q }\right) dx,
$$ 
which is negligible by applying integration by parts.
For $\frac{mN}{q^2} \ll p^{-\varepsilon}$ with some $\varepsilon>0$, we note that $\Phi\left(\frac{mN}{q^2}\right) $ is also negligible since $J_{k-1} \left(4 \pi \frac{\sqrt{xmN}}{q }\right) \ll p^{-(k-1)\varepsilon}$ and $k$ is large enough and fixed.
Hence, we may assume that $p^{-\varepsilon} \ll mN/q^2 \ll p^{\varepsilon}$. In this region, the Bessel function is not oscillatory. For $p^{-\varepsilon} \ll t \ll p^{\varepsilon}$, We can write $J_{k-1}(t)=t^{k-1}W(t)$, where $t^jW^{(j)}(t) \ll T_0$ with $T_0 \ll p^{\varepsilon}$. This is the same derivative bound satisfied by a $T_0$-inert function. By a dyadic partition for the new $m$-sum, in order to understand $\Phi\left( mN/ q^2\right)$, it suffices to consider $V \left( \frac{m}{M'}\right)$,
where $V$ is a $p^{\varepsilon}$-inert function supported on $m \asymp M'$, and $Q^2/(p^{3+\varepsilon}) \ll M' \ll Q^2/(p^{3-\varepsilon})$. 
Hence, it suffices to consider
$$
\frac{N}{q} \sum_{m}  \lambda_{f}(m) e\Big( \frac{\overline{d} m}{q} \Big) V \left( \frac{m}{M'} \right).
$$
The treatment for the $n$-sum is the same as the $m$-sum. 
For the $n$-sum
$$
\sum_{n} V\Big( \frac{n}{N}\Big)  \lambda_{f}(n)  e\Big( \frac{d n}{q} \Big),
$$
by Lemma \ref{voronoi} we obtain
$$
\frac{N}{q} \sum_{n}  \lambda_{f}(n) e\Big( \frac{-\overline{d} n}{q} \Big)  \Phi\left(\frac{n N}{q^2}\right),
$$
where  
$$
\Phi\left(\frac{nN}{q^2}\right)= 2 \pi i^k \int_0^{\infty} V(x) J_{k-1} \left(4 \pi \frac{\sqrt{xnN}}{q }\right) dx .
$$
Similarly, it suffices to consider the following
$$
\frac{N}{q} \sum_{n}  \lambda_{f}(n) e\Big( \frac{-\overline{d} n}{q} \Big) V \left( \frac{n}{N'} \right),
$$
where $V$ is a $p^{\varepsilon}$-inert function supported on $n \asymp N'$, and $Q^2/(p^{3+\varepsilon}) \ll N' \ll Q^2/(p^{3-\varepsilon})$. Hence, by applying the Voronoi summation formula to both variables, the length of the new $m,n$-sums are restricted to $Q^2/(p^{3+\varepsilon}) \ll m,n \ll Q^2/(p^{3-\varepsilon})$. By a dyadic partition, we may further assume that $N' \leq n \leq 2N'$ and $M' \leq m \leq 2M'$ with $Q^2/(p^{3+\varepsilon}) \ll M',N' \ll Q^2/(p^{3-\varepsilon})$.

 For the special case $r=0$, it is considered in Section \ref{degenerate case r=0}. By the Kuznetsov trace formula, spectral large sieve and Rankin--Selberg bound, we have the desired upper bound $p^{5+2\theta+\varepsilon}$. Next, we consider the case $r \neq 0$. Without loss of generality, we always assume that $r>0$ because of the case $r<0$ is similar. 

After applying Voronoi summation formula twice, combining with the previous result, we transfer the Ramanujan sum $S(0,n-m-r,q)$ to the Kloosterman sum $S(n-m,r,q)$.

In conclusion, it suffices to consider the following
\begin{equation}
\frac{1}{Q^2} \cdot  \sum_q V\Big(\frac{q}{Q}\Big) \cdot \frac{N^2}{q^2} \cdot \sum_{m,n} V\Big( \frac{n}{N'} \Big) V \Big( \frac{m}{M'}\Big) \lambda_{f}(n) \lambda_{f}(m)  S(n-m,r,q),
\end{equation}
where $Q^2/(p^{3+\varepsilon}) \ll M', N' \ll Q^2/(p^{3-\varepsilon})$ and $p^{3-\varepsilon} \ll N \ll p^{3+\varepsilon}$.

We assume that $n>m$ without loss of generality. 
We can further assume that $n-m > p^3 >0$. If $m=n$, we have $S(n-m,r,q)=S(0,r,q)$ and $S(0,r,q) \ll r^{1+\varepsilon}$ since $r \neq 0$ and the exact formula for the Ramanujan sum. We take the absolute value and get the upper bound $$ p^{\varepsilon} \cdot \frac{1}{Q^2} \cdot \frac{N^2}{Q^2} \cdot \frac{Q^2}{N} \cdot p^3Q \ll \frac{p^{6+\varepsilon}}{Q},$$
which is negligible since we will have $Q=p^{30}$ on the denominator. If $0<n-m \leq p^3$, by the Weil's bound for the Kloosterman sum, we have $S(n-m,r,q) \ll p^{3/2+\varepsilon}\cdot Q^{1/2}$. We take the absolute value and get the upper bound $$ p^{\varepsilon} \cdot \frac{1}{Q^2} \cdot \frac{N^2}{Q^2} \cdot \frac{Q^2}{N} \cdot p^3 \cdot p^{3/2}Q^{1/2} \cdot Q \ll \frac{p^{15/2+\varepsilon}}{Q^{1/2}},$$
which is negligible since we will have $Q^{1/2}=p^{15}$ on the denominator. 


\subsection{Kuznetsov trace formula and spectral large sieve inequality} \label{finalkuz}

We recall that it suffices to consider the following by the discussion at the end of Section \ref{twice}
\begin{equation}
\frac{1}{Q^2} \cdot  \sum_q V\Big(\frac{q}{Q}\Big) \cdot \frac{N^2}{q^2} \cdot \sum_{ n- m >p^3} V\Big( \frac{n}{N'} \Big) V \Big( \frac{m}{M'}\Big) \lambda_{f}(n) \lambda_{f}(m)  S(n-m,r,q),
\end{equation}
where $Q^2/(p^{3+\varepsilon}) \ll M', N' \ll Q^2/(p^{3-\varepsilon})$ and $p^{3-\varepsilon} \ll N \ll p^{3+\varepsilon}$.

It suffices to consider that
\begin{equation}
\frac{N^2}{Q^3} \cdot  \sum_q V\Big(\frac{q}{Q}\Big)  \sum_{ n- m >p^3} V\Big( \frac{n}{N'} \Big) V \Big( \frac{m}{M'}\Big) \lambda_{f}(n) \lambda_{f}(m) \cdot  \frac{S(n-m,r,q)}{q}.
\end{equation}

Up to a $p^{\varepsilon}$ factor, we may restrict the difference of $n$ and $m$ to $P_1 \leq n-m < 2P_1$, where $p^3 <P_1 \ll Q^2/(p^{3-\varepsilon})$. This means the following summation:
\begin{equation}
\frac{N^2}{Q^3} \cdot  \sum_q V\Big(\frac{q}{Q}\Big)  \sum_{ P_1 \leq n- m<2P_1 } V\Big( \frac{n}{N'} \Big) V \Big( \frac{m}{M'}\Big) \lambda_{f}(n) \lambda_{f}(m) \cdot  \frac{S(n-m,r,q)}{q}.
\end{equation}

We insert a $1$-inert function $W \left( \frac{\sqrt{(n-m)r}}{q} \right)$, which is supported on $0<X<x<2X$ with $X \asymp \frac{\sqrt{P_1N_4}}{Q}$. Up to a $p^{\varepsilon}$ factor, we can remove the weight function $V \left( q/Q \right)$ by a Mellin inversion (see \cite[Lemma 2.10]{Khan22}, Section \ref{subsec:inertfunctions} and \ref{weight}). Hence it suffices to consider
\begin{equation}  \label{prekuz}
\frac{N^2}{Q^3} \cdot    \sum_{ P_1 \leq n- m <2P_1} V\Big( \frac{n}{N'} \Big) V \Big( \frac{m}{M'}\Big) \lambda_{f}(n) \lambda_{f}(m) \cdot  \sum_{q \geq 1} \frac{S(n-m,r,q)}{q} \cdot W \left( \frac{\sqrt{(n-m)r}}{q} \right).
\end{equation}

Applying Kuznetsov trace formula to the inner $q$-sum (see Theorem \ref{kuznetsov} and also \cite[Theorem 2.5]{Khan22}), we can rewrite the $q$-inner sum in \eqref{prekuz} as follows: 
$$ \sum_{\substack{g \in B(1,1) \\ t_g \ll T}} \frac{1}{\vert \vert g \rvert \rvert^2} \frac{1}{\cosh({\pi t_g})} \cdot \widetilde{\phi}(t_g)\lambda_g(n-m)\overline{\lambda_g(r)}+(\cdots), $$
where $(\cdots)$ denotes the continuous spectrum contribution and the holomorphic contribution, which are similar to deal with. Since $Q^2/(p^{3+\varepsilon}) \ll m,n \ll Q^2/(p^{3-\varepsilon})$, $1 \leq r \ll p^{3+\varepsilon}$, we get $p^{-20} \ll \frac{\sqrt{P_1N_4}}{Q} \asymp \frac{\sqrt{(n-m)r}}{q} \ll p^{\varepsilon}$, whose upper bound is independent on $Q$. Applying Lemma \ref{transf-bounds} (a), we have $\widetilde{\phi}(t_g) \ll p^{\varepsilon}$ and $\widetilde{\phi}(t_g) \ll p^{-A}$ for arbitray large $A>0$ if $t_g \gg p^{\varepsilon}$. Hence, we can always restrict the spectral parameters in the Kuznetsov trace formula (see Theorem \ref{kuznetsov}) and the spectral large sieve inequality (see Theorem \ref{largesievethm}) to $k, \vert t_g \rvert , \vert t \rvert \ll p^{\varepsilon}$.

\begin{rmk}
We recall that the Selberg eigenvalue conjecture is true for level one cusp forms. Hence the contribution for the exceptional eigenvalues is void.
\end{rmk}

By the Cauchy--Schwarz inequality,  we have the upper bounds for $F_i$:
\begin{equation}\label{ub1}
    p^{\varepsilon}\cdot (N^2/Q^3) \cdot p^i \cdot \sqrt{S_1} \cdot \sqrt{S_2} \quad \quad \mathrm{for} \quad \quad i=1,2, 
\end{equation}
and 
\begin{equation}\label{ub2}
p^{\varepsilon}\cdot (N^2/Q^3) \cdot p^{i-1} \cdot \sqrt{S_1} \cdot \sqrt{S_2}  \quad \quad  \mathrm{for}  \quad \quad  i=3,
\end{equation} 
with
$$ S_1:= \sum_{\substack{g \in B(1,1) \\ t_g \ll p^{\varepsilon}}} \frac{1}{\vert \vert g \rvert \rvert^2} \frac{1}{\cosh({\pi t_g})} \cdot \left \vert \sum_{\substack{N_4/p^i \leq r_i <2N_4/p^i \\ (r_i,p)=1}} \widetilde{H}(r_i)\lambda_g(p^ir_i)\right \rvert^2+(\cdots), $$
and
$$ S_2:= \sum_{\substack{g \in B(1,1) \\ t_g \ll p^{\varepsilon}}} \frac{1}{\vert \vert g \rvert \rvert^2}\frac{1}{\cosh({\pi t_g})} \cdot \left \vert \sum_{P_1 \leq s < 2P_1} G(s)\lambda_g(s)\right \rvert^2+(\cdots),$$
where $(\cdots)$ denotes the continuous spectrum contribution and the holomorphic contribution, which are similar to deal with. Since $\lambda_j(p^i) \leq (i+1)p^{i\theta} \ll p^{i\theta+\varepsilon}$ for $i=1,2,3$ and the absolute constant is independent on the index $j$, we have $S_1 \leq (i+1)^2p^{2i\theta} \cdot \widetilde{S_1} \ll p^{2i\theta+\varepsilon} \cdot \widetilde{S_1}$, with
$$ \widetilde{S_1}:= \sum_{\substack{g \in B(1,1) \\ t_g \ll p^{\varepsilon}}} \frac{1}{\vert \vert g \rvert \rvert^2} \frac{1}{\cosh({\pi t_g})} \cdot \left \vert \sum_{\substack{N_4/p^i \leq r_i <2N_4/p^i \\ (r_i,p)=1}} \widetilde{H}(r_i)\lambda_g(r_i)\right \rvert^2+(\cdots). $$
Moreover,
$$ G(s):= \sum_{m-n=s} V\Big( \frac{n}{N'} \Big) V \Big( \frac{m}{M'}\Big) \lambda_f(m) \lambda_f(n).$$
In order to apply the large sieve inequality (Theorem \ref{largesievethm}), we need to find the upper bound for the $L^2$-norm: $\sum_{P_1 \leq s < 2P_1} \vert G(s) \rvert^2 \leq \sum_s \vert G(s) \rvert^2.$
Note that
\begin{equation}
\begin{aligned}
\sum_s \vert G(s) \rvert^2 =& \sum_{m_1-n_1=m_2-n_2} V\Big( \frac{m_1}{M'}, \frac{m_2}{M'} \Big)\times V \Big( \frac{n_1}{N'}, \frac{n_2}{N'}\Big)  \\
& \times \lambda_f(m_1)\lambda_f(n_1)\lambda_f(m_2)\lambda_f(n_2).
\end{aligned}
\end{equation}
We may replace the condition $m_1-n_1=m_2-n_2$ by a Fourier integral $\int_0^1 e(((m_1-n_1)-(m_2-n_2))\alpha)d \alpha $ and use the Wilton's bound to conclude that
$$ \sum_s \vert G(s) \rvert^2 \ll_{\varepsilon} p^{\varepsilon} \cdot M' \cdot N'.$$

Finally, by \eqref{ub1}, \eqref{ub2} and the spectral large sieve inequality (Theorem \ref{largesievethm}), we obtain  
$$ p^{1+\theta+\varepsilon} \cdot \frac{N^2}{Q^3} \cdot ((p^{\varepsilon}+p^{2+\varepsilon})\cdot p^{2+4\theta+\varepsilon})^{1/2} \cdot \left(\left(p^{\varepsilon}+\max(M',N')\right) \cdot M' \cdot N' \right)^{1/2} \ll_{\varepsilon} p^{5+\varepsilon};$$
$$ p^{2+2\theta+\varepsilon} \cdot \frac{N^2}{Q^3} \cdot ((p^{\varepsilon}+p^{1+\varepsilon})\cdot p^{1+4\theta+\varepsilon})^{1/2} \cdot \left(\left(p^{\varepsilon}+\max(M',N')\right) \cdot M' \cdot N' \right)^{1/2} \ll_{\varepsilon} p^{5+\varepsilon};$$
$$ p^{2+3\theta+\varepsilon} \cdot \frac{N^2}{Q^3} \cdot ((p^{\varepsilon}+p^{\varepsilon})\cdot p^{4\theta+\varepsilon})^{1/2} \cdot \left(\left(p^{\varepsilon}+\max(M',N')\right) \cdot M' \cdot N' \right)^{1/2} \ll_{\varepsilon} p^{5+\varepsilon},$$
for $i=1,2,3$ respectively. Combining with the estimation in Section \ref{dege}, we finish the proof for Theorem \ref{triple}. 

\end{document}